\documentclass[a4paper,11pt]{article}
\usepackage[all]{xy}
\usepackage{amsmath}
\usepackage{amsthm}
\usepackage{amsfonts}
\usepackage{amssymb}
\usepackage{amscd}
\usepackage[all]{xy}

\newcommand{\de}{\delta}
\newcommand{\lam}{\lambda}
\newcommand{\al}{\alpha}

\newcommand{\si}{\sigma}

\newcommand{\id}{\mathrm{id}}

\newcommand{\ad}{\operatorname{ad}}

\newcommand{\Hom}{\operatorname{Hom}}

\newcommand{\ot}{\otimes}

\newcommand{\trl}{\triangleleft}
\newcommand{\trr}{\triangleright}

\newcommand{\li}{{}_{1}}
\newcommand{\lii}{{}_{2}}

\newcommand{\lmo}{{}_{(0)}} 

\newcommand{\loo}{{}_{(0)}}
\newcommand{\loi}{{}_{(-1)}}

\newcommand{\lmoo}{{}_{(0)}}

\newcommand{\lmi}{{}_{(1)}}

\newcommand{\lmoi}{{}_{(-1)}}

\newcommand{\ppi}{{}_{[1]}}
\newcommand{\pii}{{}_{[2]}}

\newcommand{\qi}{{}_{<1>}}
\newcommand{\qii}{{}_{<2>}}

\def\rbiprod{{\cdot\kern-.33em\triangleright\!\!\!<}}
\def\lbiprod{{>\!\!\!\triangleleft\kern-.33em\cdot\, }}

\def\lrbiprod{{\ \cdot\kern-.60em\triangleright\kern-.33em\triangleleft\kern-.33em\cdot\, }}

\def\lprod{{>\!\!\!\triangleleft\kern-.33em\ \, }}

\newcommand{\lrcoprod}{{\,\blacktriangleright\!\!\blacktriangleleft\, }}

\newtheorem{theorem}{Theorem}[section]
\newtheorem{lemma}[theorem]{Lemma}

\newtheorem{corollary}[theorem]{Corollary}
\theoremstyle{Definition}
\newtheorem{Definition}[theorem]{Definition}
\newtheorem{remark}[theorem]{Remark}
\newtheorem{example}[theorem]{Example}

\title{Braided Hom-Lie bialgebras}
\author{Tao Zhang}
\date{}

\allowdisplaybreaks

\begin{document}
 \maketitle

 \setcounter{section}{0}

\begin{abstract}
We introduce the new concept of braided Hom-Lie bialgebras which is a generalization of Sommerh\"{a}user-Majid's braided Lie bialgebras and Yau's Hom-Lie bialgebras.  Using this concept we give the unified product construction for Hom-Lie bialgebras which can be seen as a Hom-Lie version of  Bespalov-Drabant's cocycle cross product bialgebras. Some special cases of unified products such as crossed product and matched pair of braided Hom-Lie bialgebras are investigated.
As an application, we solve the Agore-Militaru extending problem for Hom-Lie bialgebras by using some non-abelian cohomology theory.
Furthermore, one dimensional flag extending structures for Hom-Lie bialgebras are also investigated.


\end{abstract}

\tableofcontents

\section{Introduction}

In \cite{S96},  Y.~Sommerh\"{a}user introduced the concept of braided Lie bialgebras (he call it Yetter-Drinfeld Hom-Lie algebra) to give a construction of symmetrizable Kac-Moody algebras.
The theory of braided Lie bialgebras was also developed further by S.~Majid in \cite{Ma00}, where the bosonisation theorem for braided Lie bialgebras are proved.
For the theory of Yetter-Drinfel'd Hopf algebra, see \cite{BD1, BD2,R85, S02, ZC01}.
The notion of Hom-Lie bialgebras was introduced by D. Yau in \cite{yau2} which is a generalization of  Drinfeld's classical Lie bialgebra.

On the other hand, the theory of unified product and extending structure for many types of algebras were well  developed by Agore and Militaru in \cite{AM1,AM2,AM3,AM4,AM5,AM6}.
Let $A$ be a Lie (associative, Leibniz, Poisson, Jordan, etc.) algebra and $E$ a vector space containing $A$ as a subspace.
The extending problem is to describe and classify all Lie (associative, Leibniz, etc.)  algebras structures on $E$ such that $A$ is a subalgebra of $E$.
Until now, extending structures for Hom-algebra and Hom-bialgebra structures were not developed as well as the above mentioned algebra.
Recently, extending structures for 3-Lie algebras, Lie bialgebras,  infinitesimal bialgebras and Lie conformal superalgebras were studied  in \cite{zhang1,H1,zhang2,zhang3,ZCY}.

It is a  natural question whether can we develop a unified product theory  to solve the extending problem for Hom-Lie bialgebras.
In this paper,  we give a affirmation answer to this question.
The first step is to give a well definition of  braided Hom-Lie bialgebras.

The organization of  this paper is as follows.
In section \ref{sec1}, we review some basic facts and notations about Hom-Lie bialgebras and braided Hom-Lie bialgebras.
In  section  \ref{sec2},  the definition of unified product for braided Hom-Lie bialgebras is introduced.
We give the necessary and sufficient conditions for a unified product to  form Hom-Lie bialgebras.
In the last section \ref{sec3}, we study some applications of unified products, which include extending problem for braided Hom-Lie bialgebras and Hom-Lie bialgebras.
We  also study the flag extending systems.

Throughout this paper, all Hom-Lie algebras are assumed to be over an algebraically closed field $k$ of characteristic different from 2 and 3.
The space of linear maps from $(V,\al_V)$ to $W$ is denoted by $\Hom(V,W)$.
The identity map of a Hom-vector space $(V,\al_V)$ is denoted by $\id_V: V\to V$ or simply by $\id: V\to V$.
The twisting maps $\tau: V\ot V \to  V\ot V, \tau_{12}, \tau_{23}: V\ot V\ot V\to  V\ot V\ot V$ are  denoted by
$$ \tau (x\ot y)=y\ot x,\quad \tau_{12} (x\ot y\ot z)= y\ot x\ot z, \quad \tau_{23} (x\ot y\ot z)= x\ot z\ot y.$$
\section{Preliminaries}\label{sec1}

\begin{Definition} [\cite{hls}]
A Hom-Lie algebra is a triple $({H} ,[\cdot,\cdot],\alpha_H)$ where ${H}$ a vector space equipped with
a skew-symmetric brackt $[\cdot,\cdot]: {H}\otimes {H} \to {H}$  and a linear map $\alpha_H: {H} \to  {H}$ satisfying the following
Hom-Jacobi identity:
\begin{eqnarray}
&&[x,y]=-[y,x],\\
&&[\alpha_H(x),[y,z]]+[\alpha_H(y),[z,x]]+[\alpha_H(z),[x,y]]=0,
\end{eqnarray}
for all $x,y,z\in {H}$.
A Hom-Lie algebra is called a multiplicative Hom-Lie algebra if $\alpha_H$ is an algebra homomorphism, i.e, for all $x,y\in {H}$, $\alpha_H([x,y])= [\alpha_H(x), \alpha_H(y)].$
A Hom-Lie algebra is called a regular Hom-Lie algebra if $\alpha_H$ is an algebra automorphism.
\end{Definition}

By skew-symmetry of the bracket map, the Hom-Jacobi identity is equivalent to
\begin{eqnarray}
&&[\alpha_H(x),[y,z]]=[[x,y],\alpha_H(z)]+[\alpha_H(y),[x,z]].
\end{eqnarray}

For an element $x$ in a Hom-Lie algebra $({H},[\cdot, \cdot],\alpha_{H})$ and $n \geq 2$, define the adjoint map $\ad_x \colon {H}^{\otimes n} \to {H}^{\otimes n}$ by
\begin{eqnarray}
\label{eq:ad} \ad_x(y_1 \otimes \cdots \otimes y_n)
&:=& \sum_{i=1}^n \alpha_H(y_1) \otimes \cdots \otimes \alpha_H(y_{i-1}) \otimes [x,y_i] \otimes \alpha_H(y_{i+1}) \cdots \otimes \alpha_H(y_n).
\end{eqnarray}

\begin{Definition} 
A {Hom-Lie coalgebra} is a Hom-vector space $({H},\al_{H})$
equipped with a linear map $\delta: {H}\to {H} \ot {H}$, called a
cobracket, satisfying the co-anticommutativity and the co-Jacobi
identity:
\begin{enumerate}
\item[(CL1)] $\delta(x)=-\tau\delta(x)$,
\item[(CL2)] $(\al_{H} \otimes \delta)\delta(x)=(\delta\otimes \al_{H} )\delta(x)+\tau_{12}(\al_{H} \otimes \delta)\delta(x)$.
\end{enumerate}
\end{Definition}
We would like to use the sigma notation $\delta(x):=\sum x\li\ot x\lii$ for all $x\in {H}$ to denote the cobracket,
We will often omit the summation sign $\sum$ to simplify the typography.
The above conditions can be also written as:
\begin{align*}
&\mbox{(CL1)} \quad \sum  x\li\ot x\lii =- \sum  x\lii\ot x\li,\\
&\mbox{(CL2)} \quad \sum \al_H(x\li)\ot \delta(x\lii) = \sum\delta(x\li) \ot \al_H(x\lii)+\sum  \tau_{12}\left(\al_H(x\li)\ot \delta(x\lii)\right).
\end{align*}

\begin{Definition}[\cite{yau2}]\label{dfnlb}
A Hom-Lie bialgebra $(H,\al_H)$ is a Hom-vector space equipped
simultaneously with a Hom-Lie algebra structure $(H, [\cdot,\cdot], \alpha)$  and a Hom-Lie coalgebra $(H, \delta,\alpha)$ structure such that the following
compatibility condition is satisfied,
\begin{align}\label{eq:hombialg}
 \delta([x, y])=[\alpha_H(x), \delta(y)]+[\delta (x), \alpha_H(y)]=\ad_{\alpha_H(x)}\delta(y)-\ad_{\alpha_H(x)}\delta(y).
 \end{align}
 \noindent
We denote it by $(H, [\cdot,\cdot], \delta,\al_{H})$.
\end{Definition}

Using the sigma notation,  the above equation \eqref{eq:hombialg} is equivalent to
\begin{eqnarray*}
\delta([x,y])&=&\sum  [\alpha_H(x), y\li]\ot \alpha_H(y\lii)+ \sum \alpha_H(y\li)\ot[\alpha_H(x), y\lii]\\
&&+ \sum \alpha_H(x\li)\ot [x\lii,  \alpha_H(y)] + \sum [x\li,  \alpha_H(y)]\ot \alpha_H(x\lii)\\
&=&\sum  [\alpha_H(x), y\li]\ot \alpha_H(y\lii)+ \sum \alpha_H(y\li)\ot[\alpha_H(x), y\lii]\\
&&-\sum  [\alpha_H(y), x\li]\ot \alpha_H(x\lii)- \sum \alpha_H(x\li)\ot[\alpha_H(y), x\lii].
\end{eqnarray*}

A homomorphism of Hom-Lie bialgebras $\varphi: ({H}, [\cdot,\cdot],\delta)\rightarrow ({H'}, [\cdot,\cdot]',\delta')$ is both a homomorphism of Hom-Lie algebras and a homomorphism of Hom-Lie coalgebras, i.e. for all $a$, $b\in {H}$,
\begin{eqnarray*}
\varphi\circ\al_H=\al_H\circ\varphi,\quad\varphi([x,y])=[\varphi(x), \varphi(y)]',\quad \delta'\circ\varphi(x)=(\varphi\otimes \varphi)\circ \delta(x).
\end{eqnarray*}

Let $A, H$ be both Hom-Lie algebras and Hom-Lie coalgebras.   For $a, b\in
A$, $x, y\in H$,  we denote maps $$\trr:H \otimes A \to A,\quad \trl :H
\otimes A \to H, \quad \phi :A \to H \otimes A,\quad  \psi :H \to H\otimes A$$ by
\begin{eqnarray*}
&& \trr (x \otimes a) = x \triangleright a, \quad \trl(x \otimes a) = x \triangleleft a, \\
&& \phi (a)=\sum a\lmoi\ot a\loo, \quad \psi (x)= \sum x\loo \ot x\lmi.
\end{eqnarray*}

 We now fix some notations.
For a Hom-Lie algebra $(H,\al_H)$ and a linear map $\trr:H\ot A \to A$ such
that   
$$[x, y]\triangleright \alpha_A(a)=\alpha_H(x)\trr(y\trr a)-\alpha_H(y)\trr(x\trr a),$$
for all $x, y\in H, a\in A$, then $(A,
\alpha )$ is called a left $H$-Hom-Lie module. For a Hom-Lie coalgebra $(H,\al_H)$ and a
linear map $\phi :A\to H\ot A$ such that 
$$\sum \delta_H(a\loi)\ot \al_A(a\loo)=\sum  \al_H(a\loi)\ot \phi(a\loo)- \sum   \tau_{12}\left(\al_H(a\loi)\ot \phi(a\loo)\right),$$
then $(A, \phi )$ is called a left $H$-Hom-Lie comodule.
If $(H,\al_H)$ and $(A,\al_A)$ are Hom-Lie algebras, $(A,\al_A)$ is a left $H$-Hom-Lie module and
$$\al_H(x)\trr [a,b]=[x\trr a, \al_A(b)]+[\al_A(a), x\trr b],$$
then $(A, [, ], \alpha )$ is called a left $H$-Hom-Lie module algebra.
 If $(H,\al_H)$ is a Hom-Lie coalgebra and $(A,\al_A)$ is a Hom-Lie algebra,  $(A,\al_A)$ is a left $H$-Hom-Lie comodule and 
$$\phi([a, b])=\sum a\lmoi \ot [a\lmoo, b]+\sum b\lmoi\ot [a, b\lmoo],$$
then $(A,\al_A)$ is called a left $H$-Hom-Lie comodule algebra. Right Hom-Lie (co)module and Hom-Lie (co)module (co)algebra can be defined
similarly.

\begin{Definition}
Let $({A},[\cdot,\cdot])$ be a given Hom-Lie algebra (Hom-Lie coalgebra, Hom-Lie bialgebra), $(E,\al_E)$ a Hom-vector space.
An extending system of $(A,\al_A)$ through $(V,\al_V)$ is a Hom-Lie algebra  (Hom-Lie coalgebra, Hom-Lie bialgebra) on $(E,\al_E)$
such that $(V,\al_V)$ a complement subspace of $(A,\al_A)$ in $(E,\al_E)$, the canonical injection map $i: A\to E, a\mapsto (a, 0)$  or the canonical projection map $p: E\to A, (a,x)\mapsto a$ is a Hom-Lie algebra (Hom-Lie coalgebra, Hom-Lie bialgebra) homomorphism.
The extending problem is to describe and classify up to an isomorphism  the set of all Hom-Lie algebra  (Hom-Lie coalgebra, Hom-Lie bialgebra) structures that can be defined on $(E,\al_E)$.
\end{Definition}

We remark that our definition of extending system of $(A,\al_A)$ through $(V,\al_V)$ contains not only extending structure in \cite{AM1,AM2,AM3}
but also the global extension structure in \cite{AM5}.
The reason is that when we consider extending problem for Hom-Lie bialgebras, both of them are necessarily used.
Note that in our  extending system we do not demand $({A},[\cdot,\cdot])$ to be a subalgebra of $(E,\al_E)$ although  $(A,\al_A)$ is always a Hom-Lie algebra.
In fact, the canonical injection map $i: A\to E$ is a Lie (co)algebra homomorphism if and only if $(A,\al_A)$ is a Lie sub-(co)algebra of $(E,\al_E)$.

\begin{Definition}
Let $(A,\al_A)$ be a Hom-Lie algebra (Hom-Lie coalgebra, Hom-Lie bialgebra), $(E,\al_E)$ be a Hom-Lie algebra  (Hom-Lie coalgebra, Hom-Lie bialgebra) such that
${A} $ is a subspace of $(E,\al_E)$ and $(V,\al_V)$ a complement of
${A} $ in $(E,\al_E)$. Let $(E, [\cdot,\cdot])$ and $({E}, [\cdot,\cdot]')$ be two Hom-Lie algebra (Hom-Lie coalgebra, Hom-Lie bialgebra) structures on $(E,\al_E)$. For a linear map $\varphi: E \to {E}$ we consider the diagram:
\begin{equation}\label{eq:ext1}
\xymatrix{
   0  \ar[r]^{} &A \ar[d]_{\id_A} \ar[r]^{i} & E \ar[d]_{\varphi} \ar[r]^{\pi} &V \ar[d]_{\id_V} \ar[r]^{} & 0 \\
   0 \ar[r]^{} & A \ar[r]^{i'} & {E} \ar[r]^{\pi'} & V \ar[r]^{} & 0.
   }
\end{equation}
where $\pi : E \to V$ is the canonical projection of $E ={A}  \oplus V$ onto $(V,\al_V)$ and $i: {A}  \to E$ is the inclusion map. We say that $\varphi: E \to {E}$ \emph{stabilizes} $(A,\al_A)$ if the left square of the diagram \eqref{eq:ext1} is  commutative.
Two extending system $(E, [\cdot,\cdot])$ and $({E}, [\cdot,\cdot]')$ are called \emph{equivalent}, and we denote this by $(E, [\cdot,\cdot]) \equiv ({E}, [\cdot,\cdot]')$, if there exists a Hom-Lie algebra (Hom-Lie coalgebra, Hom-Lie bialgebra) isomorphism $\varphi: (E, [\cdot,\cdot])
\to ({E}, [\cdot,\cdot]')$ which stabilizes ${A} $. Denote by $Extd(E,{A} )$ ($CExtd(E,{A} )$, $BExtd(E,{A} )$) the set of equivalent classes of  Hom-Lie algebra (Hom-Lie coalgebra, Hom-Lie bialgebra) structures on $(E,\al_E)$.
\end{Definition}

\section{Braided  Hom-Lie bialgebras}
\subsection{Hom-Yetter-Drinfeld modules and braided  Hom-Lie bialgebras}

\begin{Definition}Let $(H,\al_H)$ be simultaneously a Hom-Lie alegbra and a Hom-Lie coalgebra.
 If $(V,\al_V)$ is  a left $H$-Hom-Lie module and left $H$-Hom-Lie comodule,  and satisfying the following condition:
\begin{enumerate}
\item[(YD1)] $\phi (x\trr v)= [\alpha_H(x), v_{(-1)}] \otimes \alpha\left(v_{(0)}\right)+\alpha\left(v_{(-1)}\right) \otimes \alpha_H(x) \trr v_{(0)}+\al_H(x\li) \otimes x_{2} \trr \alpha(v),$
\end{enumerate}
\noindent then $(V,\al_V)$is called a left Hom-Yetter-Drinfeld module over $(H,\al_H)$.
\end{Definition}
We denote  the  category of Hom-Yetter-Drinfeld modules over $(H,\al_H)$ by ${}^{H}_{H}\mathcal{M}$.


\begin{Definition} Let $(H,\al_H)$ be a Hom-Lie bialgebra.
If $(A,\al_A)$ be a Hom-Lie algebra, a Hom-Lie coalgebra and a left Hom-Yetter-Drinfeld module over $(H,\al_H)$, then we call $(A,\al_A)$ to be a \emph{braided Hom-Lie bialgebra} in ${}^{H}_{H}\mathcal{M}$, if the following condition is
satisfied
$$ \delta([a, b])=[\al_A(a), \delta(b)]+[\delta (a), \al_A(b)]-s(a\ot b),\leqno \mbox{(LBSa)}$$
where
\begin{eqnarray*}
s(a\ot b)&=&b_{(-1)}  \trr \al_A(a) \otimes \al_A(b_{(0)})+\al_A(a_{(0)}) \otimes  a_{(-1)}  \trr \alpha_A(b)\\
&&-a_{(-1)} \trr \alpha_A(b) \otimes \alpha_A(a_{(0)})-\al_A(b_{(0)}) \otimes  b_{(-1)}  \trr \alpha_A(a)
\end{eqnarray*}
is called the infinitesimal braiding.
\end{Definition}

When $\al=\id$,  this is exactly the Yetter-Drinfeld Lie algebra introduced  by Sommerhauser in \cite{S96} and braided Lie bialgebra studied by Majid in \cite{Ma00}.
Thus our braided Hom-Lie bialgebra is a generalization of braided Lie bialgebra.

Now we construction Hom-Lie bialgebra from braided Hom-Lie  bialgebra via biproduct.
Let $(H,\al_H)$  be an Hom-Lie bialgebra, $(A,\al_A)$ be a Hom-Lie algebra and a Hom-Lie coalgebra in ${}^{H}_{H}\mathcal{M}$.
We define bracket and cobracket on the direct sum vector space $E=A \oplus H$ by
\begin{eqnarray}
\al_E(a, x)&:=&(\al_A(a), \al_H(x)),\\
\label{eq:bracket001}
{[(a, x), (b, y)]}&:=&([a, b]+x\trr b-y\trr a, [x, y]), \\
\label{eq:cobracket001}
\delta_{E}(a, x)&:=&\delta_{A}(a)+\phi(a)-\tau \phi(a)+\delta_{H}(x).
\end{eqnarray}
This is called biproduct of $(A,\al_A)$ and $(H,\al_H)$ which will be  denoted by $A\lbiprod H$.

\begin{remark}
In sigma notation, the above equation \eqref{eq:cobracket001} can be written  as
$$\delta_{E}(a, x)= a\li\ot a\lii+a\lmoi\ot a\lmoo-a\lmoo\ot a\lmoi+x\li\ot x\lii.$$
The elements in the right hand side of this equation should be seen as in  $E\ot E=(A \oplus H)\ot (A \oplus H)\cong A\ot A \oplus H\ot A \oplus A\ot H \oplus H\ot H$,
\end{remark}

\begin{theorem}\label{boson} Let $(H,\al_H)$ is a  Hom-Lie bialgebra and
 $(A,\al_A)$  is a left Hom-Lie module and a left Hom-Lie comodule in ${}^{H}_{H}\mathcal{M}$.
Then the biproduct $A\lbiprod H$ form a Hom-Lie bialgebra if and only if  $(A,\al_A)$ is braided Hom-Lie bialgebra in ${}^{H}_{H}\mathcal{M}$.
\end{theorem}

\begin{proof} It is easy to verify that  $A\lbiprod H$ is Hom-Lie algebra and Hom-Lie coalgebra under the above bracket \eqref{eq:bracket001} and cobracket \eqref{eq:cobracket001}.
We are left to check that compatibility condition for Hom-Lie bialgebra $A\lbiprod H$:
$$
\begin{aligned}
\delta_E[(a, x), (b, y)] &=\left[\al_E(a, x),\delta_E(b, y)\right] +\left[\delta_E(a, x), \al_E(b, y)\right].
\end{aligned}
$$
The left hand side is equal to
\begin{eqnarray*}
\delta[(a, x),(b, y)]&=&\delta\left([a, b]+x \trr b-y  \trr  a ,[x, y]\right)\\
&=&\delta_{A}([a, b])+\delta_{A}\left(x\trr b\right)-\delta_{A}(y\trr  a)\\
&&+\phi([a, b])+\phi(x \trr b)-\phi\left(y  \trr  a \right)\\
&&-\tau  \phi([a, b])-\tau  \phi(x\trr b)+\tau  \phi\left(y \trr a \right)+\delta_{H}([x, y]).
\end{eqnarray*}
The right hand side is equal to
\begin{eqnarray*}
&=&\left[\alpha_A(a), b_{1}\right] \otimes \al_A(b_{2})+\alpha_H(x) \trr b_{1} \otimes \al_A(b_{2})\\
&&+\al_A(b_{1}) \otimes\left[\alpha_A(a), b_{2}\right]+\al_A(b_{1}) \otimes \alpha_H(x) \trr b_{2}\\
&&- b_{(-1)}  \trr \al_A(a) \otimes \al_A(b_{(0)})+\left[\alpha_H(x), b_{(-1)}\right] \otimes \al_A(b_{(0)})\\
&&+\al_H(b_{(-1)}) \otimes\left[\alpha_A(a), b_{(0)}\right]+\al_H(b_{(-1)}) \otimes \alpha_H(x) \trr b_{(0)}\\
&&-\left[\alpha_A(a),  b_{(0)} \right] \otimes \al_H(b_{(-1)})-\alpha_H(x) \trr  b_{(0)}  \otimes \al_H(b_{(-1)})\\
&&+\al_A(b_{(0)}) \otimes  b_{(-1)}  \trr \alpha_A(a)-\al_A(b_{(0)}) \otimes\left[\alpha_H(x),  b_{(-1)} \right]\\
&&- y_{1}  \triangleright \al_A(a) \otimes \al_H(y_{2})+\left[\alpha_H(x), y_{1}\right] \otimes \al_H(y_{2})\\
&&-\al_H(y_{1}) \otimes  y_{2} \trr \alpha_A(a)+\al_H(y_{1}) \otimes\left[\alpha_H(x), y_{2}\right]\\
&&+\left[a_{1}, \alpha_A(b)\right] \otimes \al_A(a_{2})-\alpha_H(y) \trr  a_{1}  \otimes \al_A(a_{2})\\
&&+\al_A(a_{1}) \otimes\left[a_{2}, \alpha_A(b)\right]-\al_A(a_{1}) \otimes \alpha (y)  \trr  a_{2} \\
&&+a_{(-1)}  \trr \alpha_A(b) \otimes \al_A(a_{(0)})+\left[a_{(-1)}, \alpha_H(y)\right] \otimes \al_A(a_{(0)})\\
&&+\al_H(a_{(-1)}) \otimes\left[a_{(0)}, \alpha_A(b)\right]-\al_H(a_{(-1)}) \otimes \alpha (y) \trr a_{(0)} \\
&&-\left[ a_{(0)} , \alpha_A(b)\right] \otimes \al_H(a_{(-1)})+\alpha_H(y) \trr  a_{(0)}  \otimes \al_H(a_{(-1)})\\
&&-\al_A(a_{(0)}) \otimes  a_{(-1)}  \triangleright \alpha_A(b)-\al_A(a_{(0)}) \otimes\left[ a_{(-1)} , \alpha_H(y)\right]\\
&&- x_{1}  \trr \alpha_A(b) \otimes \alpha\left(x_{2}\right)+\left[x_{1}, \alpha_H(y)\right] \otimes \alpha\left(x_{2}\right)\\
&&+\alpha\left(x_{1}\right) \otimes x_{2}  \trr \alpha_A(b)+\alpha\left(x_{1}\right) \otimes\left[x_{2}, \alpha_H(y)\right]
\end{eqnarray*}
Thus the two sides are equal to each other if and only if the following conditions hold.
\begin{eqnarray*}
(1) &&\delta_{A}\left(\left[a, b\right]\right)=\left[\alpha_A(a), b_{1}\right] \otimes \al_A(b_{2})+\al_A(b_{1}) \otimes\left[\alpha_A(a), b_{2}\right]\\
&&\qquad\qquad\qquad+\left[a_{1}, \alpha_A(b)\right] \otimes \al_A(a_{2})+\al_A(a_{1}) \otimes\left[a_{2}, \alpha_A(b)\right]\\
&&\qquad\qquad\qquad+a_{(-1)} \trr \alpha_A(b) \otimes \alpha_A(a_{(0)})+\al_A(b_{(0)}) \otimes  b_{(-1)}  \trr \alpha_A(a)\\
&&\qquad\qquad\qquad- b_{(-1)}  \trr \al_A(a) \otimes \al_A(b_{(0)})-\al_A(a_{(0)}) \otimes  a_{(-1)}  \trr \alpha_A(b),\\
(2) &&\delta_{A}\left(x\trr b\right)=\alpha_H(x) \trr b_{1} \otimes \al_A(b_{2})+\al_A(b_{1}) \otimes \alpha_H(x) \trr b_{2},\\
(3) &&\delta_{A}\left(y  \trr  a \right)=\alpha_H(y) \trr  a_{1}  \otimes \al_A(a_{2})+\al_A(a_{1}) \otimes \alpha (y) \trr a_{2} ,\\
(4) &&\phi([a, b])=\al_H(b_{(-1)}) \otimes[\alpha_A(a), b_{(0)}]+\al_H(a_{(-1)}) \otimes\left[a_{(0)}, \alpha_A(b)\right],\\
(5) &&\phi(x\trr b)=[\alpha_H(x), b_{(-1)}] \otimes \al_A(b_{(0)})+\al_H(b_{(-1)}) \otimes \alpha_H(x) \trr b_{(0)}+\al_H(x\li) \otimes x_{2} \trr \alpha_A(b),\\
(6)&&\phi\left(y \right.\left.\trr  a \right)=\left[\alpha_H(y), a_{(-1)}\right] \otimes \al_A(a_{(0)})+\al_H(a_{(-1)}) \otimes \alpha (y) \trr  a_{(0)} +\al_H(y_{1}) \otimes  y_{2} \trr \alpha_A(a),\\
(7)&&\tau  \phi([a, b])=\left[\alpha_A(a),  b_{(0)} \right] \otimes \al_H(b_{(-1)})+\left[ a_{(0)} , \alpha_A(b)\right] \otimes \al_H(a_{(-1)}),\\
(8)&&\tau  \phi(x \trr b)=\al_A(b_{(0)}) \otimes\left[\alpha_H(x),  b_{(-1)} \right]
+\alpha_H(x) \trr  b_{(0)}  \otimes \al_H(b_{(-1)})+x_{2} \trr \alpha_A(b) \otimes \alpha\left(x_{1}\right),\\
(9) &&\tau  \phi\left(y  \trr  a \right)=\al_A(a_{(0)}) \otimes\left[\alpha_H(y),  a_{(-1)} \right]+\alpha_H(y) \trr  a_{(0)}  \otimes \al_H(a_{(-1)})+ y_{1}  \trr \al_A(a) \otimes \al_H(y_{1}).
\end{eqnarray*}
From (6)--(9) we have that $A$ is a left $H$-module Hom-Lie coalgebra and $H$-comodule Hom-Lie algebra,
from (2)--(5) we get that $A$ is a left Yetter-Drinfel module over $H$, and (1) is the condition for $A$ to be a braided Hom-Lie bialgebra.
The proof is completed.
\end{proof}

\subsection{From quasitriangular Hom-Lie bialgebras to braided Hom-Lie bialgebras}
Let $(A, [\cdot,\cdot], \alpha)$ be a Hom-Lie algebra equipped with $r= r_{1}\otimes r_{2}  \in$ $A \otimes A$  such that $(\alpha \otimes \alpha)(r)=r$.
In what follows, for an element $r = \sum r_1 \otimes r_2$,  we denote  $\tau(r) = \sum r_2 \otimes r_1$.
We define the following elements in $A \otimes A \otimes A$
$$
\begin{aligned}
{[r_{12},r_{13}]}&=  [r_1,r'_1] \otimes \alpha(r_2) \otimes \alpha(r'_2), \\
{[r_{12},r_{23}]}&= \alpha\left(r_{1}\right)\otimes  [r_{2}, r'_{1}] \otimes \alpha\left(r'_{2}\right), \\
[r_{13},r_{23}]&=\alpha\left(r_{1}\right) \otimes \alpha\left(r'_{1}\right) \otimes [r'_{2}, r_{2}].
\end{aligned}
$$
where $r'=\sum_{} r'_{1} \otimes r'_{2}$ is a copy of $r$.
We consider solutions $r \in A \otimes A$ to the equation
$$
C(r)=[r_{12},r_{13}]+[r_{12},r_{23}]+[r_{13},r_{23}]=0
$$
which is called the classical Hom-Yang-Baxter equation.

\begin{Definition}\cite{yau2}
Let $(A, [\cdot,\cdot], \alpha)$ be a Hom-Lie algebra and $r \in A \otimes A$ be an element such that $(\alpha \otimes \alpha)(r)=r$. Define the linear map $\delta_{r}: A \rightarrow A \otimes A$ by
$$
\delta_{r}(x)=\left[x, r_{1}\right] \otimes \al(r_{2})+\al(r_{1})\otimes\left[x, r_{2}\right]
$$
Then
Furthermore we have $(\alpha \otimes \alpha) \circ \delta_{r}=\delta_{r} \circ \alpha$ and the following identity holds:
$$
\begin{aligned}
\delta_{r}([x,y])= [\alpha (x), \delta_{r}(y)] +  [\delta_{r}(x), \alpha (y)]
\end{aligned}
$$
Thus we obtain  $(A, [\cdot,\cdot], \delta_{r}, \alpha)$  is a Hom-Lie bialgebra.
This is called a coboundary Hom-Lie bialgebra.
A coboundary Hom-Lie bialgebra for which $C(r) =0$ is called a quasitriangular Hom-Lie bialgebra.
\end{Definition}

We remark that  $r$ is not assumed to be anti-symmetric $r+\tau(r)=0$   in a coboundary Hom-Lie bialgebra.
Note that the skew symmetry of $\delta_{r}$ is equivalent to  $r+\tau(r)$ is adjoint invariant, that is, $[x,r+\tau(r)]= 0$ for every
$x\in A$:
\begin{eqnarray}\label{eq:r-not-skew}
\left[x, r_{1}\right] \otimes \al(r_{2})+\al(r_{1})\otimes\left[x, r_{2}\right]+
\left[x, r_{2}\right]\otimes\al(r_{1})+\al(r_{2})\otimes\left[x, r_{1}\right]=0.
\end{eqnarray}

\begin{lemma}\label{lem:r01}\cite{yau2}
Let $(A,[\cdot,\cdot],\alpha,r)$ be a quasitriangular Hom-Lie bialgebra.  Then the following statements are true:
\begin{eqnarray}\label{eq:r01}
(\alpha \otimes \delta)(r) &=&[r_{13},r_{12}],\\
\label{eq:r02}
(\delta \otimes \alpha)(r) &=& [r_{13},r_{23}].
\end{eqnarray}
\end{lemma}

\begin{theorem}
Let $(A, [\cdot,\cdot], \alpha, r)$ be a quasitriangular Hom-Lie bialgebra  and $V$ be a left Hom-Lie module of $A$. Then we obtain that $V$ is a Yetter-Drinfeld module over $A$ via the left Hom-Lie comodule map $\phi: V\to A\ot V$ defined by
\begin{equation}
\phi(v)=r_2\ot r_1\trr v.
\end{equation}
\end{theorem}
\begin{proof}
First, we proof that $\phi$ is indeed a left Hom-Lie comodule:
$$\sum \delta(a\loi)\ot \al_A(a\loo)=\sum  \al (a\loi)\ot \phi(a\loo)- \sum   \tau_{12}\left(  \al (a\loi)\ot \phi(a\loo)\right).$$
Using equation \eqref{eq:r01} in the above Lemma \ref{lem:r01}, the left hand side is equal to
\begin{eqnarray*}
 \delta(a\loi)\ot \al_A(a\loo)&=&\delta_r(r_2)\ot \al(r_1)\trr \al(v)\\
&=&\alpha\left(r'_{2}\right)\ot \al(r_2)\ot [r_{1}, r'_{1}] \trr \al(v),
\end{eqnarray*}
and the hand side is equal to
\begin{eqnarray*}
&&\al (a\loi)\ot \phi(a\loo)-   \tau_{12}\left(  \al (a\loi)\ot \phi(a\loo)\right)\\
&=&\al(r_2)\ot \phi(r_1\trr v)- \tau_{12}\left(\al(r_2)\ot \phi(r_1\trr v)\right)\\
&=&\al(r_2)\ot \al(r'_2)\ot  \al(r_1)\trr (r'_1\trr v)-\al(r'_2)\ot \al(r_2)\ot  \al(r_1)\trr (r'_1\trr v)\\
&=&\al(r'_2)\ot \al(r_2)\ot \Big(\al(r_1)\trr (r'_1\trr v)- \al(r_1)\trr (r'_1\trr v)\Big)\\
&=&\alpha\left(r'_{2}\right)\ot \al(r_2)\ot [r_{1}, r'_{1}] \trr \al(v).
\end{eqnarray*}
Thus the two sides are equal to each other.

Next, we verify the Yetter-Drinfeld module condition:
\begin{eqnarray*}
\phi (x\trr v)= [\alpha (x), v_{(-1)}] \otimes \alpha\left(v_{(0)}\right)+\alpha\left(v_{(-1)}\right) \otimes \alpha (x) \trr v_{(0)}+\al (x\li) \otimes x_{2} \trr \alpha(v).
\end{eqnarray*}
Since $(\alpha \otimes \alpha)(r)=r$, we obtain that the left hand side of the above equation is equal to
\begin{eqnarray*}
\phi (x\trr v)&=&r_2\ot r_1\trr (x\trr v)=\al(r_2)\ot \al(r_1)\trr (x\trr v),
\end{eqnarray*}
and the hand side is equal to
\begin{eqnarray*}
&& [\alpha (x), v_{(-1)}] \otimes \alpha(v_{(0)})+\alpha(v_{(-1)}) \otimes \alpha (x) \trr v_{(0)}+\al (x\li) \otimes x_{2} \trr \alpha(v)\\
&=& [\alpha (x), r_2] \otimes \alpha(r_1\trr v)+\alpha(r_2) \otimes \alpha (x) \trr (r_1\trr v)-\al (x\lii) \otimes x_{1} \trr \alpha(v)\\
&=& [\alpha (x), r_2] \otimes \alpha(r_1\trr v)+\alpha(r_2) \otimes \alpha (x) \trr (r_1\trr v)\\
&&-\al(r_{2})\otimes[x, r_{1}]  \trr \alpha(v)-[\alpha (x), r_{2}]\otimes \al(r_{1}) \trr \alpha(v)\\
&=&\alpha(r_2) \otimes \alpha (x) \trr (r_1\trr v)-\al(r_{2})\otimes[x, r_{1}]  \trr \alpha(v).
\end{eqnarray*}
Therefore the two sides are equal to each other because $V$ is a Hom-Lie module over $A$: $[x, r_{1}]  \trr \alpha(v)=\alpha (x) \trr (r_1\trr v)-\al(r_1)\trr (x\trr v)$.
\end{proof}

\begin{theorem}\label{thm:transformation1}
For every quasitriangular Hom-Lie bialgebra $(A, [\cdot,\cdot], \alpha, r)$, we obtain a braided Hom-Lie bialgebra $\underline{A}$
with the new cobracket defined by
$$\underline{\delta}(x)=\al(r_{1})\otimes\left[x, r_{2}\right]+\al(r_{2})\otimes\left[x, r_{1}\right].$$
This braided Hom-Lie bialgebra is called a transmutation of the quasitriangular Hom-Lie bialgebra $(A, [\cdot,\cdot], \alpha, \delta_r)$.
\end{theorem}

\begin{proof}
We give a proof by direct computations. By the definition of $\underline{\delta}(x)$ we get
\begin{eqnarray*}
&&(\al \otimes \underline{\delta}) \underline{\delta}(x)\\
&=&(\al \otimes \underline{\delta}) (\al(r_{1}) \otimes [x, r_{2} ]+\al(r_{2}) \otimes [x, r_{1} ] )\\
&=&\al(r_{1})\otimes \al(r'_{1})\otimes [ [x, r_{2} ], r'_{2} ]+\al(r_{1})\otimes \al(r'_{2})\otimes [ [x, r_{2} ], r'_{1} ]\\
&&+\al(r_{2}) \otimes \al(r'_{1}) \otimes [ [x, r_{1} ], r'_{2} ]+\al(r_{2})\otimes \al(r'_{2})\otimes [ [x, r_{1} ], r'_{1} ]\\
&=&\al(r_{1}) \otimes \al(r'_{1})\otimes [x, [r_{2}, r'_{2} ] ]+\al(r_{1})\otimes \al(r'_{1})\otimes [ [x, r'_{2} ], r_{2} ]\\
&&+\al(r_{1}) \otimes \al(r'_{2})\otimes [x, [r_{2}, r'_{1} ] ]+\al(r_{1})\otimes \al(r'_{2})\otimes [ [x, r'_{1} ], r_{2} ]\\
&&+\al(r_{2})\otimes \al(r'_{1})\otimes [x, [r_{1}, r'_{2} ] ]+\al(r_{2})\otimes \al(r'_{1})\otimes [ [x, r'_{2} ], r_{1} ]\\
&&+\al(r_{2})\otimes \al(r'_{2})\otimes [x, [r_{1}, r'_{1} ] ]+\al(r_{2})\otimes \al(r'_{2})\otimes [ [x, r'_{1} ], r_{1} ]\\
&=&-\al(r_{1}) \otimes \al(r'_{2})\otimes [x, [r_{2}, r'_{1} ] ]-\al(r_{1}) \otimes [r_{2}, r'_{1} ] \otimes [x, r'_{2} ]\\
&&-\al(r_{1})\otimes [r_{2}, r'_{2} ] \otimes [x, r'_{1} ]+\al(r_{1})\otimes \al(r'_{1})\otimes [ [x, r'_{2} ], r_{2} ]\\
&&+\al(r_{1})\otimes \al(r'_{2})\otimes [x, [r_{2}, r'_{1} ] ]+\al(r_{1})\otimes \al(r'_{2})\otimes [ [x, r'_{1} ], r_{2} ]\\
&&-\al(r_{2})\otimes \al(r'_{2})\otimes [x, [r_{1}, r'_{1} ] ]-\al(r_{2})\otimes [r_{1}, r'_{1} ] \otimes [x, r'_{2} ]\\
&&-\al(r_{2})\otimes [r_{1}, r'_{2} ] \otimes [x, r'_{1} ]+\al(r_{2})\otimes \al(r'_{1})\otimes [ [x, r'_{2} ], r_{1} ]\\
&&+\al(r_{2})\otimes \al(r'_{2})\otimes [x, [r_{1}, r'_{1} ] ]+\al(r_{2})\otimes \al(r'_{2})\otimes [ [x, r'_{1} ], r_{1} ]\\
&=&-\al(r_{1}) \otimes [r_{2}, r'_{1} ] \otimes [x, r'_{2} ]-\al(r_{1})\otimes [r_{2}, r'_{2} ] \otimes [x, r'_{1} ]\\
&&+\al(r_{1})\otimes \al(r'_{1})\otimes [ [x, r'_{2} ], r_{2} ]+\al(r_{1})\otimes \al(r'_{2})\otimes [ [x, r'_{1} ], r_{2} ]\\
&&-\al(r_{2})\otimes [r_{1}, r'_{1} ] \otimes [x, r'_{2} ]-\al(r_{2})\otimes [r_{1}, r'_{2} ] \otimes [x, r'_{1} ]\\
&&+\al(r_{2})\otimes \al(r'_{1})\otimes [ [x, r'_{2} ], r_{1} ]+\al(r_{2})\otimes \al(r'_{2})\otimes [ [x, r'_{1} ], r_{1} ]
\end{eqnarray*}
where we use the Hom-Jacobi identity in the third equality and \eqref{eq:r-not-skew} in the fourth equality.
In the meantime, we have
\begin{eqnarray*}
&&\tau_{12}(\al \otimes \underline{\delta}) \underline{\delta}(x)+(\underline{\delta}(\otimes\id) \underline{\delta}(x)\\
&=&\tau_{12}(\al \otimes \underline{\delta}) (\al(r_{1})  \otimes [x, r_{2} ]+\al(r_{2})  \otimes [x, r_{1} ] )\\
&&+(\underline{\delta}\otimes \al) (\al(r_{1}) \otimes [x, r_{2} ]+\al(r_{2})  \otimes [x, r_{1} ] )\\
&=&\al(r'_{1})\otimes \al(r_{1})\otimes [ [x, r_{2} ], r'_{2} ]+\al(r'_{2})\otimes \al(r_{1})\otimes [ [x, r_{2} ], r'_{1} ] \\
&&+\al(r'_{1})\otimes \al(r_{2})\otimes [ [x, r_{1} ], r'_{2} ]+\al(r'_{2})\otimes \al(r_{2})\otimes [ [x, r_{1} ], r'_{1} ] \\
&&+\al(r'_{1})\otimes [r_{1}, r'_{2} ] \otimes [x, r_{2} ]+\al(r'_{2})\otimes [r_{1}, r'_{1} ] \otimes [x_{1}, r_{2} ] \\
&&+\al(r'_{1})\otimes [r_{2}, r'_{2} ] \otimes [x, r_{1} ]+\al(r'_{2})\otimes [r_{2}, y_{1}^{1} ] \otimes [x, r_{1} ] \\
&=&\al(r_{1})\otimes \al(r'_{1})\otimes [ [x, r'_{2} ], r_{2} ]+\al(r_{2})\otimes \al(r'_{1})\otimes [ [x, r'_{2} ], r_{1} ] \\
&&+\al(r_{1})\otimes \al(r'_{2})\otimes [[x, r'_{1}], r_{2} ]+\al(r_{2})\otimes \al(r'_{2})\otimes [ [x, r'_{1} ], r_{1} ] \\
&&+\al(r_{1})\otimes [r'_{1}, r_{2} ] \otimes [x, r'_{2} ]+\al(r_{2})\otimes [r'_{1}, r_{1} ] \otimes [x, r'_{2} ] \\
&&+\al(r_{1})\otimes [r'_{2}, r_{2} ] \otimes [x, r'_{1} ]+\al(r_{2})\otimes [r'_{2}, r_{1} ] \otimes [x, r'_{1} ]
\end{eqnarray*}
Thus the two sides are equal to each other by comparing term by term using skew-symmetry of Hom-Lie algebra $A$.
\end{proof}

\begin{example}
Let $A=\mathfrak{sl}(2)=span\{H, X, Y\}$ be the three dimensional Hom-Lie algebra with
$\al: \mathfrak{sl}(2)\to \mathfrak{sl}(2)$ given by
$$\al(H)=H,\quad \al(X)=kX,\quad\al(Y)=k^{-1}Y,$$
the Lie bracket  given by
\begin{eqnarray*}
[H, X]=2kX,\quad [H,Y]=-2k^{-1}Y,\quad [X,Y]=H
\end{eqnarray*}
and the Lie cobracket given by
\begin{eqnarray*}
\delta(H)=0,~~~\delta(X)=kX\wedge H,~~~~\delta(Y)=k^{-1}Y\wedge H.
\end{eqnarray*}
Then $(\mathfrak{sl}(2),\al)$ is a quasi-triangular Hom-Lie bialgebra with $r=X\ot Y+ \frac{1}{4}H\ot H.$
By the above Theorem \ref{thm:transformation1}, we obtain a braided Hom-Lie bialgebra with the new cobracket by
\begin{eqnarray*}
\underline{\delta}(H)=2k^2(X\otimes Y+Y\otimes X),~~~\underline{\delta}(X)=kX\wedge H,~~~~\underline{\delta}(Y)=k^{-1}Y\wedge H.
\end{eqnarray*}
\end{example}

\section{Unified product for braided Hom-Lie bialgebras}\label{sec2}

In this section, we will construct unified product for braided Hom-Lie bialgebras.
First, we review the notion of matched pairs of Hom-Lie algebras and Hom-Lie coalgebras.

\subsection{Unified product of braided Hom-Lie bialgebras}
In the following definitions, we introduce some new concept of cocycle Hom-Lie algebras and cycle Hom-Lie coalgebras, which are  in fact not really ordinary Hom-Lie algebras and Hom-Lie coalgebras, but generalized ones.
Using these new algebraic structures, we will construct the unified product of braided Hom-Lie bialgebras.
Denote maps
$$\si: H\ot H\to A,\quad \theta: A\ot A\to H,\quad P: A\to H\ot H,\quad Q: H\to A\ot A$$
 by
\begin{align*}
&\si(x, y)\in A,\quad\theta(a, b)\in H,\\
&P(a)=\sum a\ppi\ot a\pii\in H\ot H,\quad Q(x)=\sum x\qi\ot x\qii\in A\ot A.
 \end{align*}

An antisymmetric bilinear map $\si: H\ot H\to A$ is called a cocycle on $(H,\al_H)$ if
$$\al_H(x)\trr\si(y,z)+\al_H(y)\trr\si(z,x)+\al_H(z)\trr\si(x,y)=\si([x, y],\al_H(z))+\si([y, z], \al_H(x))+\si([z, x], \al_H(y)). \leqno\mbox{(CC1)}$$
An antisymmetric bilinear map $\theta: A\ot A\to H$ is called a cocycle on $(A,\al_A)$ if
$$\theta(a, b)\trl \al_A(c)+\theta(b, c)\trl \al_A(a)+\theta(c,a)\trl \al_A(b)=\theta(\al_A(a), [b, c])+\theta(\al_A(b), [c, a])+\theta(\al_A(c), [a,b]).\leqno\mbox{(CC2)}$$
A co-antisymmetric linear map $P: A\to H\ot H$ is called a cycle on $(A,\al_A)$ if
\begin{align*}
\mbox{(CC3)}& \quad  \al_H(a\lmoi)\ot P(a\lmoo)+\tau_{12}\tau_{23}\left(\al_H(a\lmoi)\ot P(a\lmoo)\right)+\tau_{23}\tau_{12}\left(\al_H(a\lmoi)\ot P(a\lmoo)\right)\qquad\quad\qquad\quad\\
&\qquad  =\delta(a\ppi)\ot \al_H(a\pii)+\tau_{12}\tau_{23}\left(\delta(a\ppi)\ot \al_H(a\pii)\right)+\tau_{23}\tau_{12}\left(\delta(a\ppi)\ot \al_H(a\pii)\right).
\end{align*}
A co-antisymmetric linear map $Q: H\to A\ot A$ is called a cycle on $(H,\al_H)$ if
\begin{align*}
\mbox{(CC4)}& \quad  Q(x\lmoo)\ot \al_A(x\lmi) +\tau_{12}\tau_{23}\left(Q(x\lmoo)\ot \al_A(x\lmi) \right)+\tau_{23}\tau_{12}\left(Q(x\lmoo)\ot \al_A(x\lmi) \right)\\
&= \al_A(x\qi)\ot \delta(x\qii)+\tau_{12}\tau_{23}\left( \al_A(x\qi)\ot \delta(x\qii) \right)+\tau_{12}\tau_{23}\left(\al_A(x\qi)\ot \delta(x\qii)\right).\quad
\end{align*}

\begin{Definition}\label{def1234}
(i): Let $\si:H\ot H\to A$ be a cocycle on  $(H,\al_H)$ equipped with an  antisymmetric bilinear map  $[\cdot,\cdot]: H \ot H \to H$, satisfying the the following cocycle Jacobi identity:
$$[[x, y], \al_H(z)]+[[y, z], \al_H(x)]+[[z, x], \al_H(y)]= \al_H(x)\trl \si(y, z)+\al_H(y)\trl \si(z, x)+\al_H(z)\trl \si(x, y).\leqno\mbox{(CC5)}$$
Then  $(H,\al_H)$ is called a  $\si$-Hom-Lie algebra.

(ii): Let $\theta: A\ot A\to H$ be a cocycle on $(A,\al_A)$  equipped with  antisymmetric bilinear map  $[\cdot,\cdot]: A \ot A \to A$, satisfying the the following cocycle Jacobi identity:
$$ [[a, b], \al_A(c)]+[[b, c], \al_A(a)]+[[c, a], \al_A(b)]= \theta(a,b)\trr \al_A(c)+\theta(b, c)\trr \al_A(a)+\theta(c, a)\trr \al_A(b).\leqno\mbox{(CC6)}$$
\noindent Then  $(A,\al_A)$ is called a  $\theta$-Hom-Lie algebra.

(iii) Let $P: A\to H\ot H$ be a cycle on  $(H,\al_H)$ equipped with a co-antisymmetric linear map $\delta: H \to H \ot H$, satisfying the the
following cycle co-Jacobi identity:
\begin{align*}
\mbox{(CC7)}& \quad \delta(x\li)\ot \al_H(x\lii)+\tau_{12}\tau_{23}\left(\delta(x\li)\ot \al_H(x\lii) \right)+\tau_{23}\tau_{12}\left( \delta(x\li)\ot \al_H(x\lii)\right)\\
& \quad \qquad = \al_H(x\lmoo)\ot P(x\lmi)+\tau_{12}\tau_{23}\left(\al_H(x\lmoo)\ot P(x\lmi)  \right)+\tau_{23}\tau_{12}\left(\al_H(x\lmoo)\ot P(x\lmi) \right).\qquad\quad
\end{align*}
\noindent Then  $(H,\al_H)$ is called a  $P$-Hom-Lie coalgebra.

(iv) Let $Q: H\to A\ot A$ be a cycle on  $(A,\al_A)$ equipped with a co-antisymmetric linear map $\delta: A \to A \ot A$, satisfying the the
following cycle co-Jacobi identity:
\begin{align*}
\mbox{(CC8)} &\quad \delta(a\li)\ot \al_A(a\lii)+\tau_{12}\tau_{23}\left(\delta(a\li)\ot\al_A(a\lii)\right)+\tau_{23}\tau_{12}\left( \delta(a\li)\ot \al_A(a\lii)\right)\qquad\qquad\\
& \quad =  Q(a\lmoi)\ot \al_A(a\lmoo)+\tau_{12}\tau_{23}\left( Q(a\lmoi)\ot \al_A(a\lmoo)\right)+\tau_{23}\tau_{12}\left( Q(a\lmoi)\ot \al_A(a\lmoo) \right).
\end{align*}
\noindent Then  $(A,\al_A)$ is called a $Q$-Hom-Lie coalgebra.
\end{Definition}

First we consider the Hom-Lie algebra structures on $E=A\oplus H$.

\begin{theorem} \label{lem001}
Let $(A,\al_A)$ be a $\theta$-Hom-Lie algebra and $(H,\al_H)$ be a  $\sigma$-Hom-Lie algebra.
Then  $E=A\oplus H$ is a Hom-Lie algebra with bracket given by:
\begin{align}
\al_E(a, x)&=(\al_A(a), \al_H(x)),\\
[(a, x), (b, y)]&=([a, b]+x\trr b-y\trr a+\sigma(x, y), [x,y]+x\trl b-y\trl a+\theta(a, b))
\end{align}
if and only if the following compatibility conditions hold:
\begin{enumerate}
\item[(TM1)] $[x, y]\trr \al_A(a)+ [\si(x, y), \al_A(a)]=\al_H(x)\trr (y\trr a)-\al_H(y)\trr (x\trr a)+ \si(\al_H(x), y\trl a)+\si(x\trl a, \al_H(y)),$
\item[(TM2)] $\al_H(x)\trl [a, b]+[\al_H(x), \theta(a, b)]=(x\trl a)\trl \al_A(b)-(x\trl b)\trl \al_A(a) +\theta(x\trr a, \al_A(b))+\theta(\al_A(a), x\trr b),$
\item[(TBB1)] $ \al_H(x)\trr [a, b]+\si(\al_H(x), \theta(a, b))=[x\trr a, \al_A(b)]+[\al_A(a), x\trr b]+(x\trl a)\trr \al_A(b)- (x\trl b)\trr \al_A(a),$
\item[(TBB2)] $[x, y]\trl \al_A(a)+\theta(\si(x, y), \al_A(a))=[\al_H(x), y\trl a]+[x\trl a, \al_H(y)]+\al_H(x)\trl (y\trr a)- \al_H(y)\trl(x\trr a).$
\end{enumerate}
In this case, $(A,H)$ is called a cocycle cross product system. This Hom-Lie algebra will be denoted by $A_{\trr, \theta}\# {}_{\trl, \sigma}H$
\end{theorem}

\begin{proof}
We have to check when
\begin{equation}\label{eq:cocycleproduct001}
[\al_H(x), [a, b]_{E}]_{E}+[\al_A(a), [b, x]_{E}]_{E}+ [\al_A(b), [x, a]_{E}]_{E}=0.
\end{equation}
In fact,
\begin{eqnarray*}
&&[\al_H(x), [a,b]_{E}]_{E}=\al_H(x)\trr [a, b]+\al_H(x)\trl [a, b]+[\al_H(x), \theta(a, b)]+\si(\al_H(x), \theta(a, b)),\\
&&[\al_A(a), [b,x]_{E}]_{E}=-[\al_A(a), x\trr b]-\theta(\al_A(a), x\trr b)+(x\trl b)\trr \al_A(a)+ (x\trl b)\trl \al_A(a), \\
&&[\al_A(b), [x, a]_{E}]_{E}=[\al_A(b), x\trr a]+\theta(\al_A(b), x\trr a)-(x\trl a)\trr \al_A(b)-(x\trl a) \trl \al_A(b).
 \end{eqnarray*}
 Thus equation \eqref{eq:cocycleproduct001} holds if and only if  (TM2) and (TBB1) hold.
 Similarly, one can verify that $[\al_A(a), [x, y]_{E}]_{E}+[\al_H(x), [y, a]_{E}]_{E}+ [\al_H(y), [a, x]_{E}]_{E}=0$ holds if and only if (TM1) and (TBB2) hold.
\end{proof}

Next we consider the Hom-Lie coalgebra structure on $E=A\oplus H$.
\begin{theorem} \label{lem3}
Let $(A,\al_A)$ be a $Q$-Hom-Lie coalgebra and $(H,\al_H)$ be  a $P$-Hom-Lie coalgebra. If we define $E=A^{\phi, P}\# {}^{\psi, Q} H$ as the vector
space $A\oplus H$ with the Lie cobracket
\begin{eqnarray}\label{eqcobracket01}
\delta_{E}(a)=\delta_{A}(a)+\phi(a)-\tau\phi(a)+P(a), \quad\delta_{E}(x)=\delta_{H}(x)+\psi(x)-\tau\psi(x)+Q(x),
\end{eqnarray}
then  $A^{\phi, P}\# {}^{\psi, Q} H$ is a Hom-Lie coalgebra if and only if the following compatibility conditions hold:
\begin{enumerate}
\item[(TM3)] $ \delta_H(a\lmoi)\ot \al_A(a\lmoo)+P(a\li)\ot  \al_A(a\lii)$\\
$= \al_H(a\lmoi)\ot \phi(a\lmoo) - \tau_{12}\left( \al_H(a\lmoi)\ot \phi(a\lmoo)\right)$\\
$ + \al_H(a\ppi)\ot \psi\left(a\pii\right)  +\tau_{23}\left(\psi(a\ppi)\ot \al_H(a\pii)\right)$,
\item[(TM4)] $  \al_H(x\lmoo) \ot \delta_A(x\lmi)+\al_H(x\li)\ot Q(x\lii)$\\
$= \psi(x\lmoo)\ot \al_A(x\lmi) -  \tau_{23}\left(\psi(x\lmoo)\ot \al_A(x\lmi)\right)$\\
$+  \psi(x\qi)\ot \al_A(x\qii) + \tau_{12}\left(\al_A(x\qi)\ot\psi( x\qii)\right)$,
\item[(TBB3)] $\al_H(a\loi)\ot \delta_A(a\loo) +\al_H(a\ppi)\ot Q(a\pii)$\\
$=  \phi(a\li)\ot \al_A(a\lii)+  \tau_{12}\left( \al_A(a\li)\ot \phi(a\lii)\right)$\\
$+  \psi(a\loi)\ot \al_A(a\loo)-  \tau_{23}\left( \psi(a\loi)\ot \al_A(a\loo)\right)$,
\item[(TBB4)] $\delta_H(x\lmo)\ot \al_A(x\lmi)+P(x\qi)\ot \al_A(x\qii)$\\
$= \al_H(x\li)\ot \psi(x\lii) +  \tau_{23}\left( \psi(x\li)\ot \al_H(x\lii)\right)$\\
$+  \al_H(x\lmo)\ot \psi(x\lmi) -  \tau_{12}\left( \al_H(x\lmo)\ot \psi(x\lmi)\right)$.
 \end{enumerate}
 In this case, $(A,H)$ is called a cycle cross coproduct system.
\end{theorem}

The proof Theorem \ref{lem3} is dual to Theorem \ref{lem001} so we omit the details.

Combing the two Theorem \ref{lem001}  is dual to Theorem \ref{lem3} together, we  obtain an ordinary Hom-Lie bialgebra from two cocylce braided Hom-Lie bialgebras.
\begin{theorem} \label{main2}
Let $(A, H)$ be a cocycle cross product system and a cycle cross
coproduct system. Then the Hom-Lie algebra  $A_{\trr, \theta}\# {}_{\trl, \sigma}H$ and
Hom-Lie coalgebra $A^{\phi, P}\# {}^{\psi, Q} H$ fit together to form an ordinary
Hom-Lie bialgebra if and only if the following compatibility conditions hold:
\begin{enumerate}
\item[(TBB5)] $\delta_{A} (x\trr a)+Q(x\trl a)=\al_H(x)\trr a\li\ot \al_A(a\lii)+\al_A(a\li)\ot \al_H(x)\trr a\lii$\\
$\qquad+ x\loo\trr \al_A(a)\ot \al_A(x\lmi)-\al_A(x\lmi)\ot x\loo\trr\al_A(a)+[Q(x), \al_A(a)]$\\
$\qquad +\si(\al_H(x), a\lmoi)\ot \al_A(a\lmoo)- \al_A(a\lmoo)\ot\si(\al_H(x), a\lmoi)$,
\item[(TBB6)] $ \delta_{H} (x\trl a)+P(x\trr a)=\al_H(x\li)\ot x\lii \trl \al_A(a)+x\li\trl \al_A(a)\ot x\lii$\\
$+ \al_H(a\loi)\ot \al_H(x)\trl a\loo-\al_H(x)\trl a\loo\ot \al_H(a\loi)$\\
$+[\al_H(x), P(a)]+\al_H(x\lmoo)\ot \theta(x\lmi, \al_A(a))-\theta(x\lmi, \al_A(a))\ot \al_H(x\lmoo)$,
\item[(TBB7)] $  \phi([a, b])+\psi\theta(a, b)=\al_H(a\loi)\ot [a\loo, \al_A(b)]+\al_H(b\loi)\ot [\al_A(a), b\loo]$\\
$+a\loi\trl \al_A(b) \ot \al_A(a\loo) -b\loi\trl \al_A(a) \ot \al_A(b\loo)$\\
$+\theta(\al_A(a), b\li)\ot \al_A(b\lii)+\theta(a\li, \al_A(b))\ot \al_A(a\lii)$\\
$+\al_H(a\ppi)\ot a\pii\trr \al_A(b)-\al_H(b\ppi)\ot b\pii\trr \al_A(a)$,
\item[(TBB8)] $  \psi([x, y])+\phi\si([x, y])=[\al_H(x), y\loo]\ot \al_A(y\lmi)+[x\loo, \al_H(y)]\ot \al_A(x\lmi)$\\
$+\al(y\loo)\ot \al_H(x)\trr y\lmi -\al_H(x\loo)\ot \al_H(y)\trr x\lmi$\\
$+ \al_H(x)\li\ot\si(x\lii, \al_H(y))+\al_H(y\li)\ot \si(\al_H(x), y\lii)$\\
$+\al_H(x)\trl y\qi\ot \al_A(y\qii)-\al_H(y)\trl x\qi\ot \al_A(x\qii)$,
\item[(TLB1)] $\delta_H\theta(a, b)+P([a,b])= \al_H(a\lmoi)\ot\theta(a\lmoo, \al_A(b))+\al_H(b\lmoi)\ot\theta(\al_A(a),b\lmoo)$\\
$- \theta(\al_A(a), b\lmoo)\ot \al_H(b\lmoi) -\theta(a\lmoo, \al_A(b))\ot \al_H(a\lmoi)$\\
$+\al_H(a\ppi)\ot a\pii \trl \al_A(b)+a\ppi\trl \al_A(b) \ot \al_H(a\pii)$\\
$-b\ppi\trl \al_A(a)\ot \al_H(b\pii) -\al_H(b\ppi)\ot b\pii\trl \al_A(a)$,
\item[(TLB2)] $\delta_A \si(x, y)+Q([x, y])=\si(x\lmoo,  y)\ot \al_A(x\lmi)+\si(\al_H(x), y\lmoo)\ot \al_A(y\lmi)$\\
$-\al_A(y\lmi)\ot \si(\al_H(x), y\lmoo)-\al_A(x\lmi)\ot \si(x\lmoo, y)$\\
$+\al_H(x)\trr y\qi\ot \al_A(y\qii)+\al_A(y\qi)\ot \al_H(x)\trr y\qii $\\
$-\al_A(x\qi)\ot  \al_H(y)\trr  x\qii - \al_H(y)\trr x\qi \ot \al_A(x\qii)$,
\item[(TLB3)] $\delta_A([a, b])+Q\theta(a, b)=[\al_A(a), \delta_A(b)] +[\delta_A(a), \al_A(b)]$\\
$-b_{(-1)}  \trr \al_A(a) \otimes \al_A(b_{(0)})-\al_A(a_{(0)}) \otimes  a_{(-1)}  \trr \alpha_A(b)$\\
$+a_{(-1)} \trr \alpha_A(b) \otimes \alpha_A(a_{(0)})+\al_A(b_{(0)}) \otimes  b_{(-1)}  \trr \alpha_A(a),$
\item[(TLB4)] $\delta_H([x, y])+P\sigma(x, y) = [\delta_H(x), \al_H(y)]+[\al_H(x), \delta_H([y)] $\\
$-y_{(-1)}  \trr \alpha (x) \otimes \al_H(y_{(0)})-\al_H(x_{(0)}) \otimes  x_{(-1)}  \trr \alpha_H(y)$\\
$+x_{(-1)} \trr \alpha_H(y) \otimes \al_H(x_{(0)})+\al_H(y_{(0)}) \otimes  y_{(-1)}  \trr \alpha_H(x)$,
\item[(TYD)] $\phi (x\trr a)+\psi (x\trl a)=$\\
 $[\alpha_H(x), a_{(-1)}] \otimes \al_A(a_{(0)})+\al_H(a_{(-1)}) \otimes \alpha_H(x) \trr a_{(0)}+\al_H(x\li) \otimes x_{2} \trr \alpha_A(a)$\\
 $+\al_H(x\loo)\ot [x\lmi, \al_A(a)]+x\loo\trl \al_A(a)\ot \al_A(x\lmi) +\al_H(x)\trl a\li\ot \al_A(a\lii)$\\
 $+\al_H(a\ppi)\ot \si(\al_H(x), a\pii)+\al_A(x\qi)\ot\theta(x\qii, \al_A(a))$.
\end{enumerate}
This Hom-Lie bialgebra is denote by $A^{\phi, P}_{\trr, \theta}\# {}^{\psi, Q}_{\trl, \sigma}H$.
We call it the unified product of $(A,\al_A)$ and $(H,\al_H)$.
\end{theorem}

\begin{proof}
We investigate the case  of (LB) on   $A
\otimes A$:
\begin{eqnarray*}
&&\delta_{E}([a, b]_{E})=\delta_{E}([a, b]+\theta(a, b))=\\
&=&\delta([a, b])+\phi([a, b])-\tau\phi([a, b])+P([a, b])\\
&&+\delta\theta(a, b)+\psi\theta (a, b)-\tau\psi\theta(a,
b)+Q\theta(a, b).
\end{eqnarray*}
Denote the right hand side terms by $(a), (b), \cdots , (h).$
\begin{eqnarray*}
 && [\al_A(a), \delta_{E}(b)] + [\delta_{E}(a), \al_A(b)]\\
&=& [\al_A(a),  (\delta+\phi-\tau\phi+P)(b)] + [(\delta+\phi-\tau\phi+P)(a), \al_A(b)]\\
&=&[\al_A(a), b\li]\ot \al_A(b\lii)(1)+\al_A(b\li)\ot[\al_A(a), b\lii](2) -b\lmoi\trr \al_A(a)\ot \al_A(b\lmoo)(3)\\
&&-b\lmoi\trl \al_A(a)\ot \al_A(b\lmoo)(4) +\al_H(b\lmoi)\ot [\al_A(a), b\lmoo](5)- [\al_A(a), b\lmoo]\ot \al_H(b\lmoi)(6)\\
&& +\al_A(b\lmoo)\ot b\lmoi\trr \al_A(a)(7)+ \al_A(b\lmoo)\ot b\lmoi\trl \al_A(a)(8)\\
&& +\al_H(b\lmoi)\ot\theta(\al_A(a), b\lmoo)(19)-\theta(\al_A(a), b\lmoo)\ot \al_H(b\lmoi)(20)-b\ppi\trr \al_A(a)\ot \al_H(b\pii)(21)\\
&&-b\ppi\trl \al_A(a)\ot \al_H(b\pii)(22) -\al_H(b\ppi)\ot b\pii\trr \al_A(a)(23)-\al_H(b\ppi)\ot b\pii\trl \al_A(a)(24)\\
&&+\al_A(a\li)\ot [a\lii, \al_A(b)](9)+[a\li, \al_A(b)]\ot \al_A(a\lii)(10) \\
&&+\al_H(a\lmoi)\ot [a\lmoo, \al_A(b)](11)+ a\lmoi \trr \al_A(b)\ot \al_A(a\lmoo)(12)+a\lmoi\trl \al_A(b) \ot \al_A(a\lmoo)(13)\\
&&-\al_A(a\lmoo)\ot a\lmoi\trr \al_A(b)(14)- \al_A(a\lmoo)\ot a\lmoi\trl \al_A(b) (15)\\
&&- [a\lmoo, \al_A(b)]\ot \al_H(a\lmoi)(16)+\theta(\al_A(a), b\li)\ot \al_A(b\lii)(17)+ \al_A(b\li)\ot \theta(\al_A(a), b\lii)(18)\\
&&+\al_A(a\li)\ot \theta(a\lii,\al_A(b))(25)+ \theta(a\li,\al_A(b))\ot \al_A(a\lii)(26)+\al_H(b\lmoi)\ot\theta(a\lmoo, \al_A(b))(27)\\
&&-\theta(a\lmoo,\al_A(b))\ot \al_H(b\lmoi)(28)+\al_H(a\ppi)\ot  a\pii \trr \al_A(b)(29)+\al_H(a\ppi)\ot a\pii \trl \al_A(b)(30) \\
&&+ a\ppi\trr \al_A(b)\ot \al_H(a\pii) (31)+a\ppi\trl \al_A(b) \ot \al_H(a\pii)(32).
\end{eqnarray*}

\noindent Then by (TLB3) we get $(a)+(h)=(1)+(2)
+(9)+(10)-(3)-(14)+(12)+(7)$; by (TBB7) we get
$(b)+(f)=(11)+(5)+(13)-(4)+(17)+(26)+(29)-(23)$,
$(c)+(g)=(16)+(6)+(8)-(13)-(18)-(25)-(31)+(21)$; by (TLB1) we
get$(d)+(e)=(27)+(19)-(20)-(28)+(30)+(32)-(22)-(24)$.

We investigate the case  of (LB) on   $H \otimes H$:
\begin{eqnarray*}
&&\delta_{E}([x, y]_{E})=\delta_{E}([x, y]+\si(x, y))=\\
&=&\delta([x, y])+\psi([x, y])-\tau\psi([x, y])+Q([x, y])\\
&&+\delta\si(x, y)+\phi\si (x, y)-\tau\phi\si(x, y)+P\si(x, y)
\end{eqnarray*}
Denote the right hand side terms by $(a), (b), \cdots , (h)$.
\begin{eqnarray*}
 && [\al_H(x), \de_{E}(y)]+ [\de_{E}(x), \al_H(y)]\\
&=& [\al_H(x), (\de+\psi-\tau\psi+Q)(y)]+[ (\de+\psi-\tau\psi+Q)(x), \al_H(y)]\\
&=&[\al_H(x), y\li]\ot \al_H(y\lii)(1)+\al_H(y\li)\ot[\al_H(x), y\lii](2)\\
&&+[\al_H(x), y\lmoo]\ot \al_A(y\lmi)(3)+\al_H(y\lmoo)\ot \al_H(x)\trr y\lmi(4)+\al_H(y\lmoo)\ot \al_H(x)\trl y\lmi(5)\\
&& -\al_H(x)\trr y\lmi\ot \al_H(y\lmoo)(6)-\al_H(x)\trl y\lmi \ot \al_H(y\lmoo)(7)- \al_A(y\lmi)\ot [\al_H(x), y\lmoo](8)\\
&&+\si(\al_H(x), y\li)\ot \al_H(y\lii)(17)+ \al_H(y\li)\ot \si(\al_H(x), y\lii)(18)\\
&&+\si(\al_H(x), y\lmoo)\ot \al_A(y\lmi)(19)-\al_A(y\lmi)\ot \si(\al_H(x), y\lmoo)(20)\\
&&+\al_H(x)\trr y\qi\ot \al_A(y\qii)(21)+\al_H(x)\trl y\qi\ot \al_H(y)\qii(22)\\
&&+\al_A(y\qi)\ot \al_H(x)\trr y\qii(23) +\al_A(y\qi)\ot \al_H(x)\trl y\qii(24)\\
&&+\al_H(x\li)\ot [x\lii, \al_H(y)](9)+[x\li, \al_H(y)]\ot \al_H(x\lii)(10)\\
&&-\al_H(x\lmoo)\ot \al_H(y)\trr x\lmi(11)- \al_H(x\lmoo)\ot \al_H(y)\trl x\lmi(12)+[x\lmoo, \al_H(y)]\ot \al_A(x\lmi)(13)\\
&&-\al_A(x\lmi)\ot [x\lmoo, \al_H(y)](14)+\al_H(y)\trr x\lmi  \ot \al_H(x\lmoo)(15)+\al_H(y)\trl  x\lmi \ot \al_H(x\lmoo)(16)\\
&&+\al_H(x\li)\ot \si(x\lii, \al_H(y))(25)+ \si(x\li, \al_H(y))\ot \al_H(x\lii)(26)\\
&&+\si(x\lmoo,  y)\ot \al_A(x\lmi)(27)-\al_A(x\lmi)\ot \si(x\lmoo, y)(28)\\
&&-\al_A(x\qi)\ot  \al_H(y)\trr  x\qii(29)-\al_A(x\qi)\ot  \al_H(y)\trl  x\qii (30)\\
&& - \al_H(y)\trr x\qi \ot \al_A(x\qii)(31) -\al_H(y)\trl  x\qi \ot \al_A(x\qii)(32)
\end{eqnarray*}
\noindent Then by (TLB4) we get $(a)+(h)=(1)+(2)
+(9)+(10)-(12)-(7)+(5)+(16)$; by (TBB8) we get
$(b)+(f)=(3)+(13)+(4)-(11)+(25)+(18)+(22)-(32)$,
$(c)+(g)=(8)+(14)+(6)-(15)-(26)-(17)-(24)+(30)$; by (TLB2) we
get$(d)+(e)=(27)+(19)-(20)-(28)+(21)+(23)-(29)-(31)$.

We now check the axiom  (LB) on $H\otimes A$. For $x\in H, a\in A$,
we get the equality below:
\begin{eqnarray*}
&&\delta_{E}([x, a]_{E}) =\delta_{E}(x\trr a)+\delta_{E}(x\trl a)=\\
&=&\delta_{A}(x\trr a)+\phi(x\trr
a)-\tau\phi(x\trr a)+P(x\trr a)\\
&&+\delta_{H}(x\trl a)+\psi(x\trl a)-\tau\psi(x\trl a)+Q(x\trl a)
\end{eqnarray*}
Denote the right hand side terms by $(a), (b), \cdots , (h)$.
\begin{eqnarray*}
&&[\al_H(x), \delta_{E}(a)]+ [\delta_{E}(x), \al_A(a)]\\
&=& [\al_H(x), (\delta+\phi-\tau\phi+P)(a) +  (\delta+\psi-\tau\psi+Q)(x), \al_A(a)]\\
&=& [\al_H(x),  a\li]\ot \al_A(a\lii)(1) +\al_H(x)\trl a\li\ot \al_A(a\lii) (2)+\al_A(a\li)\ot \al_H(x)\trr a\lii(3)\\
&&+\al_A(a\li)\ot \al_H(x)\trl a\lii(4)+[\al_H(x), a\lmoi]\ot \al_A(a\lmoo)(5)+\al_H(a\lmoi)\ot \al_H(x)\trr a\lmoo(6)\\
&&+ \al_H(a\lmoi)\ot \al_H(x)\trl a\lmoo(7)-\al_H(x)\trr a\lmoo \ot \al_H(a\lmoi)(8)-\al_H(x)\trl a\lmoo\ot \al_H(a\lmoi)(9)\\
&&- \al_A(a\lmoo)\ot [\al_H(x), a\lmoi](10)\\
&&+\si(\al_H(x), a\lmoi)\ot \al_A(a\lmoo)(21)- \al_A(a\lmoo)\ot\si(\al_H(x), a\lmoi)(22)\\
&&+[\al_H(x), a\ppi]\ot \al_H(a\pii)(23)+\si(\al_H(x), a\ppi)\ot \al_H(a\pii)(24)\\
&&+\al_H(a\ppi)\ot [\al_H(x), a\pii](25)+\al_H(a\ppi)\ot \si(\al_H(x), a\pii)(26)\\
&&+\al_H(x\li)\ot x\lii\trr \al_A(a) (11)+\al_H(x\li)\ot x\lii\trl \al_A(a)(12)\\
&&+x\li\trr \al_A(a)  \ot \al_H(x\lii) (13)+x\li\trl a\ot \al_H(x\lii)(14)\\
&&+ \al_H(x\lmoo)\ot [x\lmi, \al_A(a)]  (15)+ x\lmoo\trr \al_A(a)\ot \al_A(x\lmi) (16)+x\lmoo\trl a\ot \al_A(x\lmi)(17)\\
&&-\al_A(x\lmi)\ot x\lmoo\trr \al_A(a)  (18) - \al_A(x\lmi) \ot x\lmoo\trl \al_A(a)(19)- [x\lmi, \al_A(a)]\ot \al_H(x\lmoo)(20)\\
&&+\al_H(x\lmoo)\ot\theta(x\lmi, \al_A(a))(27)-\theta(x\lmi, \al_A(a))\ot \al_H(x\lmoo)(28)\\
&&+\al_A(x\qi)\ot \theta(x\qii, \al_A(a))(29)+\al_A(x\qi)\ot [x\qii, \al_A(a)](30)\\
&&+[x\qi, \al_A(a)]\ot \al_A(x\qii)(31)+\theta(x\qi, \al_A(a))\ot \al_A(x\qii) (32)
\end{eqnarray*}
\noindent Then by (TBB5) we get $(a)+(h)=(1)+(3)
+(16)-(18)+(29)+(31)+(21)-(22)$;  by (TBB6) we
get$(d)+(e)=(12)+(14)+(7)-(8)+(23)+(25)+(27)-(28)$; by (TYB) we get
$(b)+(f)=(5)+(6)+(11)+(15)+(17)+(2)+(26)+(32)$,
$(c)+(g)=(10)+(9)-(13)+(20)+(19)-(4)-(26)-(30)$.
\end{proof}

In the case $\theta=0,  P=0$, then $(A, [\cdot,\cdot])$ is a Hom-Lie algebra and $(A, \delta_A)$ is a Hom-Lie coalgebra and by (TLB3) we obtain that $(A,[\cdot,\cdot], \delta_A)$ is a braided Hom-Lie bialgebra in ${}^{H}_{H}\mathcal{M}$.
In the case $\sigma=0,  Q=0$, then $(H, [\cdot,\cdot], \trr)$ is a Hom-Lie algebra and $(H, \delta_H)$ is a Hom-Lie coalgebra and by (TLB4) we obtain that $(H,[\cdot,\cdot], \delta_H)$ is  really a braided Hom-Lie bialgebra in ${\mathcal{M}}^{A}_{A}$.
That is why in Theorem \ref{main2} we call  $A^{\phi, P}_{\trr, \theta}\# {}^{\psi, Q}_{\trl, \sigma}H$ the unified product for braided Hom-Lie bialgebras.

Put $\theta=0, Q=0$, then from (TLB3) we get that $(A,\al_A)$ is a  braided Hom-Lie bialgebra. By the above Theorem \ref{main2}, we obtain:
\begin{theorem}\label{thm1}
Let $(A,\al_A)$ be a braided Hom-Lie bialgebra and $(V,\al_V)$ a Hom-vector space.
An extending datum of $(A,\al_A)$ by $(V,\al_V)$ is  $\Omega^b({A},V)=(\trr, \trl, \sigma, [\cdot,\cdot], \phi, \psi, Q, \delta_V)$ consisting of  eight linear maps
\begin{eqnarray*}
\trr: V\times {A}\rightarrow {A},~~~~\trl: V\times {A}\rightarrow V,~~~~\sigma:  V\times V \rightarrow {A},~~~[\cdot,\cdot]:V\times V \rightarrow V,\\
 \phi :A \to V\otimes A, \quad{\psi}: V\to  V\otimes A,~~~~{P}: A\rightarrow {V}\otimes {V},~~~~\delta_V: V\rightarrow V\otimes V.
\end{eqnarray*}
Then the unified product $A^{\phi}_{\trr}\# {}^{\psi,Q}_{\trl, \sigma}\, V$ with bracket
\begin{align}
[(a, x), (b, y)]:=([a, b]+x\trr b-y\trr a+\sigma(x, y), [x,y]+x\trl b-y\trl a)
\end{align}
and cobracket
\begin{eqnarray}
\delta_E(a)=\delta_{A}(a)+{\phi}(a)-\tau{\phi}(a)+P(a),\quad \delta_E(x)=\delta_V(x)+{\psi}(x)-\tau{\psi}(x)
\end{eqnarray}
form a Hom-Lie bialgebra if and only if $A{}_{\,\trr}\# {}_{\trl, \sigma} V$ form a Hom-Lie algebra, $A^{\phi}\# {}^{\psi,Q} \, V$ form a Hom-Lie coalgebra and the following conditions are satisfied:
\begin{enumerate}
\item[(TBB5)] $\delta_{A} (x\trr a)=\al_H(x)\trr a\li\ot \al_A(a\lii)+\al_A(a\li)\ot \al_H(x)\trr a\lii$\\
$\qquad+ x\loo\trr \al_A(a)\ot \al_A(x\lmi)-\al_A(x\lmi)\ot x\loo\trr\al_A(a)$\\
$\qquad +\si(\al_H(x), a\lmoi)\ot \al_A(a\lmoo)- \al_A(a\lmoo)\ot\si(\al_H(x), a\lmoi)$,
\item[(TBB6)] $ \delta_{H} (x\trl a)+P(x\trr a)=\al_H(x\li)\ot x\lii \trl \al_A(a)+x\li\trl \al_A(a)\ot x\lii$\\
$+ \al_H(a\loi)\ot \al_H(x)\trl a\loo-\al_H(x)\trl a\loo\ot \al_H(a\loi)+[\al_H(x), P(a)]$,
\item[(TBB7)] $  \phi([a, b])=\al_H(a\loi)\ot [a\loo, \al_A(b)]+\al_H(b\loi)\ot [\al_A(a), b\loo]$\\
$+a\loi\trl \al_A(b) \ot \al_A(a\loo) -b\loi\trl \al_A(a) \ot \al_A(b\loo)$\\
$+\al_H(a\ppi)\ot a\pii\trr \al_A(b)-\al_H(b\ppi)\ot b\pii\trr \al_A(a)$,
\item[(TBB8)] $  \psi([x, y])+\phi\sigma([x,y])=[\al_H(x), y\loo]\ot \al_A(y\lmi)+[x\loo, \al_H(y)]\ot \al_A(x\lmi)$\\
$+\al(y\loo)\ot \al_H(x)\trr y\lmi -\al_H(x\loo)\ot \al_H(y)\trr x\lmi$\\
$+ \al_H(x)\li\ot\si(x\lii, \al_H(y))+\al_H(y\li)\ot \si(\al_H(x), y\lii)$,
\item[(TLB1)] $P([a,b])= \al_H(a\ppi)\ot a\pii \trl \al_A(b)+a\ppi\trl \al_A(b) \ot \al_H(a\pii)$\\
$-b\ppi\trl \al_A(a)\ot \al_H(b\pii) -\al_H(b\ppi)\ot b\pii\trl \al_A(a)$,
\item[(TLB2)] $\delta_A \si(x, y)=\si(x\lmoo,  y)\ot \al_A(x\lmi)+\si(\al_H(x), y\lmoo)\ot \al_A(y\lmi)$\\
$-\al_A(y\lmi)\ot \si(\al_H(x), y\lmoo)-\al_A(x\lmi)\ot \si(x\lmoo, y)$,
\item[(TLB3)] $\delta_A([a, b])=[\al_A(a), \delta_A(b)] +[\delta_A(a), \al_A(b)]$\\
$-b_{(-1)}  \trr \al_A(a) \otimes \al_A(b_{(0)})-\al_A(a_{(0)}) \otimes  a_{(-1)}  \trr \alpha_A(b)$\\
$+a_{(-1)} \trr \alpha_A(b) \otimes \alpha_A(a_{(0)})+\al_A(b_{(0)}) \otimes  b_{(-1)}  \trr \alpha_A(a),$
\item[(TLB4)] $\delta_H([x, y])+P\sigma(x, y) = [\delta_H(x), \al_H(y)]+[\al_H(x), \delta_H([y)] $\\
$-y_{(-1)}  \trr \alpha (x) \otimes \al_H(y_{(0)})-\al_H(x_{(0)}) \otimes  x_{(-1)}  \trr \alpha_H(y)$\\
$+x_{(-1)} \trr \alpha_H(y) \otimes \al_H(x_{(0)})+\al_H(y_{(0)}) \otimes  y_{(-1)}  \trr \alpha_H(x)$,
\item[(TYD)] $\phi (x\trr a)+\psi (x\trl a)=$\\
 $[\alpha_H(x), a_{(-1)}] \otimes \al_A(a_{(0)})+\al_H(a_{(-1)}) \otimes \alpha_H(x) \trr a_{(0)}+\al_H(x\li) \otimes x_{2} \trr \alpha_A(a)$\\
 $+\al_H(x\loo)\ot [x\lmi, \al_A(a)]+x\loo\trl \al_A(a)\ot \al_A(x\lmi) +\al_H(x)\trl a\li\ot \al_A(a\lii)$\\
 $+\al_H(a\ppi)\ot \si(\al_H(x), a\pii)$.
\end{enumerate}
\end{theorem}

\subsection{Matched pairs of braided Hom-Lie bialgebras}
In the following of this section, we investigate  the spacial case when the cocycle maps $\sigma,\theta$ are zero. In this case, we obtain Hom-Lie bialgebras from two braided Hom-Lie bialgebras. The other cases for which $(A,\al_A)$ is a Hom-Lie bialgebra will be given in the next sections.

\begin{Definition}
Assume that  $(A,\al)$ and  $(H,\al)$ are Hom-Lie algebras. If   $(A, \trr )$ is a left $H$-Hom-Lie module, $(H, \trl)$ is a right $A$-module, and the following condition (M1) and (M2) hold, then
 $(A, H, \trr, \trl )$ ( or  $(A, H )$) is called  a\emph{ matched pair of  Lie  algebras}:
\begin{enumerate}
\item[(BB1)] $\al_H(x)\trr [a, b]=[x\trr a, \al_A(b)]+[ \al_A(a), x\trr b]+(x\trl a)\trr \al_A(b)- (x\trl b)\trr \al_A(a),$
\item[(BB2)] $[x, y]\trl \al_A(a)=[\al_H(x), y\trl a]+[x\trl a, \al_H(y)]+\al_H(x)\trl (y\trr a)- \al_H(y)\trl(x\trr a).$
\end{enumerate}
\end{Definition}

\begin{lemma}\label{lem1}
Let $(A, H)$ be a matched pair of Hom-Lie algebras,  then
we obtain a new Hom-Lie algebra on the vector space $E=A\oplus H$ with
bracket given by
\begin{align*}
\al_E(a, x)&=(\al_A(a), \al_H(x)),\\
[(a, x), (b, y)]&=([a, b]+x\trr b-y\trr a, [x, y]+x\trl b-y\trl a).
\end{align*}
We will denote it by $A\bowtie H$. \end{lemma}

The dual version is the matched pair of Hom-Lie coalgebras.
\begin{Definition}
Two Hom-Lie coalgebras $( A, H)$ form a \emph{matched pair of Hom-Lie coalgebras} if $(A, \phi)$ is a left $H$-Hom-Lie comodule and $(H, \psi)$ is
a right $A$-comodule, obeying the conditions
\begin{enumerate}
\item[(BB3)] $ (\al \ot \delta)\phi= ( \phi \ot \al)\delta +(\tau\ot \al) ( \al\ot \phi)\delta + ( \psi \ot \al)\phi + (\al\ot\tau) ( \psi \ot \al)\phi,$
\item[(BB4)] $ (\delta \ot \al)\psi= ( \al \ot \psi)\delta +(\al\ot \tau) ( \al\ot \psi)\delta + ( \al \ot \phi)\psi +  (\tau\ot\al)( \al\ot \phi)\psi.$
\end{enumerate}
In sigma notation, the above conditions are
\begin{align*}
\mbox{(BB3)} \quad \sum \al_H(a\loi)\ot \delta_A(a\loo) =&\sum \phi(a\li)\ot \al_A(a\lii)+ \sum\tau_{12}\left( \al_A(a\li)\ot \phi(a\lii)\right)\\
&+\sum \psi(a\loi)\ot \al_A(a\loo)-\sum \tau_{23}\left( \psi(a\loi)\ot \al_A(a\loo)\right),
 \end{align*}
 \begin{align*}
\mbox{(BB4)} \quad  \sum \delta_H(x\lmo)\ot \al_A(x\lmi)=&\sum \al_H(x\li)\ot \psi(x\lii) +\sum \tau_{23}\left( \psi(x\li)\ot \al_H(x\lii)\right)\\
&+\sum \al_H(x\lmo)\ot \psi(x\lmi) -\sum \tau_{12}\left( \al_H(x\lmo)\ot \psi(x\lmi)\right).\quad
 \end{align*}
\end{Definition}

\begin{lemma} \label{lem2}
 Let $(A, H)$ be a matched pair of Hom-Lie coalgebras. We define $E=A\lrcoprod H$ as the vector space $A\oplus H$ with Lie cobracket
$$\delta_{E}(a)=(\delta_{A}+\phi-\tau\phi)(a),\quad\delta_{E}(x)=(\delta_{H}+\psi-\tau\psi)(x),$$
that is
$$\delta_{E}(a)=\sum a\li \ot a\lii+\sum a\loi \ot a\loo-\sum a\loo\ot a\loi, $$
$$\delta_{E}(x)=\sum x\li \ot x\lii+\sum x\loo \ot x\lmi-\sum  x\lmi \ot x\loo.$$
Then  $A\lrcoprod H$ is a Hom-Lie coalgebra.
\end{lemma}

\begin{Definition}
Let $(A, H)$ be matched pair of Hom-Lie algebras and Hom-Lie coalgebras.
If  the following conditions hold:
\begin{enumerate}
\item[(BB5)] $\delta_{A} (x\trr a)=\al_H(x)\trr a\li\ot \al_A(a\lii)+\al_A(a\li)\ot \al_H(x)\trr a\lii$\\
$\qquad+ x\loo\trr \al_A(a)\ot \al_A(x\lmi)-\al_A(x\lmi)\ot x\loo\trr\al_A(a),$
\item[(BB6)] $ \delta_{H} (x\trl a)=\al_H(x\li)\ot x\lii \trl \al_A(a)+x\li\trl \al_A(a)\ot x\lii$\\
$\qquad+ \al_H(a\loi)\ot \al_H(x)\trl a\loo-\al_H(x)\trl a\loo\ot \al_H(a\loi),$
\item[(BB7)] $  \phi([a, b])=\al_H(a\loi)\ot [a\loo, \al_A(b)]+\al_H(b\loi)\ot [\al_A(a), b\loo]$\\
$+a\loi\trl \al_A(b) \ot \al_A(a\loo) -b\loi\trl \al_A(a) \ot \al_A(b\loo),$
\item[(BB8)] $  \psi([x, y])=[\al_H(x), y\loo]\ot \al_A(y\lmi)+[x\loo, \al_H(y)]\ot \al_A(x\lmi)$\\
$+\al(y\loo)\ot \al_H(x)\trr y\lmi -\al_H(x\loo)\ot \al_H(y)\trr x\lmi,$
\item[(BB9)] $\delta_A([a, b])=[\al_A(a), \delta_A(b)] +[\delta_A(a), \al_A(b)]$\\
$-b_{(-1)}  \trr \al_A(a) \otimes \al_A(b_{(0)})-\al_A(a_{(0)}) \otimes  a_{(-1)}  \trr \alpha_A(b)$\\
$+a_{(-1)} \trr \alpha_A(b) \otimes \alpha_A(a_{(0)})+\al_A(b_{(0)}) \otimes  b_{(-1)}  \trr \alpha_A(a),$
\item[(BB10)] $\delta_H([x, y])= [\delta_H(x), \al_H(y)]+[\al_H(x), \delta_H([y)] -x\lmoo\ot \al_H(y)\trl x\lmi -\al_H(x)\trl y\lmi\ot y\lmoo$\\
$+ y\lmoo\ot \al_H(x)\trl y\lmi+\al_H(y)\trl x\lmi\ot x\lmoo$,
\item[(YDB)] $ \phi (x\trr a)+\psi (x\trl a)
=  [\alpha_H(x), a_{(-1)}] \otimes \alpha_A(a_{(0)})+\alpha_H(a_{(-1)}) \otimes \alpha_H(x) \trr a_{(0)}+\al_H(x\li) \otimes x_{2} \trr \alpha_A(a)$\\
$+\al_H(x\loo)\ot [x\lmi, \al_A(a)]+x\loo\trl \al_A(a)\ot \al_A(x\lmi) +\al_H(x)\trl a\li\ot \al_A(a\lii),$
\end{enumerate}
then $(A, H)$ is called a \emph{double matched pair}.
\end{Definition}

\begin{theorem}\label{main1}
Let $(A, H)$ be matched pair of Hom-Lie algebras and Hom-Lie coalgebras.
If we define the double biproduct of $(A,\al_A)$ and
$(H,\al_H)$, denoted by $A\lrbiprod H$, $A\lrbiprod H=A\bowtie H$ as Lie
algebra, $A\lrbiprod H=A\lrcoprod H$ as Hom-Lie coalgebra, then
$A\lrbiprod H$ become a Hom-Lie bialgebra if and only if $(A, H)$ is a double matched pair.
 \end{theorem}

 In particular, the condition (YDB) can be splitting into
 \begin{enumerate}
\item[(YDB1)] $\phi (x\trr a)=  [\al_H(x), a_{(-1)}] \ot \al(a_{(0)})+\al_H(a_{(-1)}) \ot\al_H(x) \trr a_{(0)}+\al_H(x_1) \ot x_{2} \trr \al_A(a)$
    \end{enumerate}
 and
 \begin{enumerate}
\item[(YDB2)] $ \psi (x\trl a)=\al_H(x\loo)\ot [x\lmi, \al_A(a)]+x\loo\trl \al_A(a)\ot \al_A(x\lmi) +\al_H(x)\trl a\li\ot \al_A(a\lii).$
\end{enumerate}
In this case, $(A,\al_A)$  is a left Hom-Yetter-Drinfeld module in ${}^{H}_{H}\mathcal{M}$ and $(H,\al_H)$  is a right Hom-Yetter-Drinfeld module in $\mathcal{M}^{A}_{A}$.
Together with (BB8) and (BB9), we obtain that $(A,\al_A)$ is  a  braided Hom-Lie bialgebra in ${}^{H}_{H}\mathcal{M}$, $(H,\al_H)$ is  a braided Hom-Lie bialgebra in $\mathcal{M}^{A}_{A}$.

\section{Applications}\label{sec3}

In this section, we will study the extending problem and non-abelian extension problem for  Hom-Lie bialgebra.
We will find some special cases when the braided Hom-Lie bialgebra $({A}, [\cdot,\cdot], \delta_{A})$ is reduced to an ordinary Hom-Lie bialgebra.
It is proved that these problems can be solved by using the non-abelian cohomology theory based on our unified product for braided Hom-Lie bialgebras in last section. The proof of most result of this section is by direct computation, so we omit them. For more details, the reader could see \cite{zhang2}.

\subsection{Extending structures for Hom-Lie algebras}

There are two cases for $(A,\al_A)$  to be a Hom-Lie algebra. The first case is when $\trr=0, \theta\neq 0$,  from (TBB1) we get $\si(x, \theta(a, b))=0$, since
$\theta\neq 0$ we assume $\sigma=0$ for simplicity, thus  we obtain the following type (a1)  unified product for Hom-Lie algebras.
\begin{corollary}\label{cor01}
Let $({A},[\cdot,\cdot])$ be a Hom-Lie algebra and $(V,\al_V)$ a Hom-vector space. An extending datum of $(A,\al_A)$ by $(V,\al_V)$ of type (a1)  is  $\Omega^{(1)}({A},V)=(\beta, \theta, [\cdot,\cdot]_V)$ consisting of bilinear maps
\begin{eqnarray*}
\theta: A\times {A}\rightarrow {V},~~~~\trl: V\times {A}\rightarrow V,~~~[\cdot,\cdot]_V:V\times V \rightarrow V.
\end{eqnarray*}
Denote by $A_{\theta}\#_{\trl}V$ the vector space $E={A}\oplus V$ with  bracket $[\cdot,\cdot]: E\times E \rightarrow E$ given by
\begin{eqnarray}
[(a, x), (b, y)]:=\left([a,b],\,  x\trl b-y\trl a+[x,y]+\theta(a,b)\right),
\end{eqnarray}
for all $a$, $b\in {A}$, $x$, $y\in V$.
Then $A_{\theta}\# {}_{\trl}V$ is a Hom-Lie algebra if and only if the following compatibility conditions hold for all $a$, $b\in {A}$, $x$, $y$, $z\in V$:
\begin{enumerate}
\item[(A1)] $\theta(a,a)=0,~~~{[x,x]}=0,$
\item[(A2)] $[[x, y], \al_V(z)]+[[y, z], \al_V(x)]+[[z, x], \al_V(y)]=0,$
\item[(A3)] $[x, y]\trl \al_A(a)=[\al_V(x), y\trl a]+[x\trl a, \al_V(y)],$
\item[(A4)] $\al_V(x)\trl [a, b]+[\al_V(x), \theta(a, b)]=(x\trl a)\trl \al_A(b)-(x\trl b)\trl \al_A(a),$
\item[(A5)] $\theta(a, b)\trl \al_A(c)+\theta(b, c)\trl \al_A(a)+\theta(c,a)\trl \al_A(b)=\theta(\al_A(a), [b, c])+\theta(\al_A(b), [c, a])+\theta(\al_A(c), [a,b]).$
\end{enumerate}
\end{corollary}
Note that in this case by (A2)  we obtain that $(V,\al_V)$ is a Hom-Lie algebra. Furthermore, $(V,\al_V)$ is in fact a subalgebra of $A_{\theta}\#_{\trl} V$  but $(A,\al_A)$ is not. Instead $(A,\al_A)$ is a quotient algebra of $A_{\theta}\#_{\trl} V$.

Denote the set of all Hom-Lie algebra extending datum of $(A,\al_A)$ by $(V,\al_V)$ of type (a1)  by $\mathcal{A}^{(1)}({A},V)$.

Note that $A_{\theta}\#_{\trl} V$  is a Hom-Lie algebra containing ${V}$ as a subalgebra.
In fact any Hom-Lie algebraic structure on $(E,\al_E)$ containing $(A,\al_A)$  as subspace and ${V}$ as subalgebra is isomorphic to such a unified product of this type.

In the following, we always assume that $(A,\al_A)$ is a subspace of a Hom-vector space $(E,\al_E)$, there exists a projection map $p: E \to{A}$ such that $p(a) = a$, for all $a \in {A}$.
Then the kernel space $V := \ker(p)$ is also a subspace of $(E,\al_E)$ and a complement of $(A,\al_A)$ in $(E,\al_E)$.

\begin{lemma}\label{lem:33-1}
Let $({A},[\cdot,\cdot])$ be a Hom-Lie algebra and $(E,\al_E)$ a Hom-vector space containing $(A,\al_A)$ as a subspace.
Suppose that there is a Hom-Lie algebraic structure $(E, [\cdot,\cdot]_E)$ on $(E,\al_E)$ such that $(V,\al_V)$ is a Lie subalgebra of $(E,\al_E)$
and the canonical projection map $p: E\to A$ is a Hom-Lie algebra homomorphism.
Then there exists a Hom-Lie algebraic extending datum $\Omega^{(1)}({A},V)$ of $(A,\al_A)$ by $(V,\al_V)$ such that
$(E, [\cdot,\cdot]_E)\cong A_{\theta}\#_{\trl} V$.
\end{lemma}

\begin{proof}
Since $(V,\al_V)$ is a subalgebra of $(E,\al_E)$, we have $[x,y]_E\in V$.
We define the extending datum of $(A,\al_A)$ through $(V,\al_V)$ by the following formulas:
\begin{eqnarray*}
\trl: V\times {A} \to V, \qquad {x} \triangleleft {a} &:=&[{x}, {a}]_E,\\
\theta: A\times A \to V, \qquad \theta(a,b) &:=&[a,b]_E-p \bigl([a,b]_E\bigl),\\
{[\cdot,\cdot]_V}: V \times V \to V, \qquad [{x}, {y}]_V&:=& [{x}, {y}]_E .
\end{eqnarray*}
for any $a , b\in {A}$ and $x, y\in V$. It is easy to see that the above maps are  well defined and
$\Omega^{(1)}({A}, V) = \bigl(\theta, \trl,[\cdot,\cdot]_V\bigl)$ is an extending system of
${A}$ through $(V,\al_V)$ and
\begin{eqnarray*}
\varphi:A_{\theta}\#_{\trl} V\to E, \qquad \varphi(a, x) := a+x
\end{eqnarray*}
is an isomorphism of Hom-Lie algebras.
\end{proof}

\begin{lemma}
Let $\Omega^{(1)}({A},V)=(\theta, \trl, [\cdot,\cdot]_V)$ and ${\Omega}'^{(1)}({A},V)=(\theta', \trl', [\cdot,\cdot]_V')$ be two Hom-Lie algebraic extending datums of $(A,\al_A)$ by $(V,\al_V)$ of type (a1) and $A_{\theta}\#_{\trl} V$, $A_{\theta'}\#_{\trl'} V$ be the corresponding unified products. Then there exists a bijection between the set of all homomorphisms of Hom-Lie algebras $\varphi:A_{\theta}\#_{\trl} V\to A_{\theta'}\#_{\trl'} V$ whose restriction on $(A,\al_A)$ is the identity map and the set of pairs $(r,s)$, where $r:V\rightarrow {A}$ and $s:V\rightarrow V$ are two linear maps satisfying
\begin{eqnarray}
&&{r}(x\trl a)=[{r}(x),a],\\
&&[a,b]'=[a,b]+r\theta(a,b),\\
&&{r}([x, y])=[{r}(x),{r}(y)]',\\
&&{s}(x)\trl' a+\theta'(r(x), a)={s}(x\trl a),\\
&&\theta'(a,b)=s\theta(a,b),\\
&&{s}([x, y])=[{s}(x),{s}(y)]'+{s}(x)\trl'{r}(y)-{s}(y)\trl'{r}(x)+\theta'(r(x), r(y)),
\end{eqnarray}
for all $a, b\in{A}$ and $x$, $y\in V$.

Under the above bijection the homomorphism of Hom-Lie algebras $\varphi=\varphi_{r,s}: A_{\theta}\#_{\trl} V\to A_{\theta'}\#_{\trl'} V$ to $(r,s)$ is given  by $\varphi(a,x)=(a+r(x), s(x))$ for all $a\in {A}$ and $x\in V$. Moreover, $\varphi=\varphi_{r,s}$ is an isomorphism if and only if $s: V\rightarrow V$ is a linear isomorphism.
\end{lemma}

The second case is when $\theta=0, \trr\neq 0$,  we obtain the following type (a2)  unified product for Hom-Lie algebras.

\begin{corollary}\label{cor02}
Let $(A,\al_A)$ be a Hom-Lie algebra and $(V,\al_V)$ a Hom-vector space. An extending datum of $(A,\al_A)$ by $(V,\al_V)$ of  type (a2)  is  $\Omega^{(2)}({A},V)=(\trr, \trl, \sigma, [\cdot,\cdot])$ consisting of four bilinear maps
\begin{eqnarray*}
\trr: V\times {A}\rightarrow {A},~~~~\trl: V\times {A}\rightarrow V,~~~~\sigma:  V\times V \rightarrow {A},~~~[\cdot,\cdot]:V\times V \rightarrow V.
\end{eqnarray*}
Denote by $A{}_{\,\trr}\# {}_{\trl, \sigma}H$ the vector space $E={A}\oplus V$ with the bilinear map
$[\cdot,\cdot]: E\times E \rightarrow E$ given by
\begin{eqnarray}
[(a, x), (b, y)]:=\left([a,b]+x\trr b-y\trr a+\sigma(x,y),\,  x\trl b-y\trl a+[x,y]\right),
\end{eqnarray}
for all $a$, $b\in {A}$, $x$, $y\in V$.
Then $A{}_{\,\trr}\# {}_{\trl, \sigma}V$ is a Hom-Lie algebra if and only if the following compatibility conditions hold for all $a$, $b\in {A}$, $x$, $y$, $z\in V$:
\begin{enumerate}
\item[(B1)] $\sigma(x,x)=0,~~~{[x,x]}=0,$
\item[(B2)] $\al_V(x)\trl [a,b]=(x\trl a)\trl \al_A(b)-(x\trl b)\trl \al_A(a),$
\item[(B3)] $\al_V(x)\trr [a, b]=[x\trr a, \al_A(b)]+[ \al_A(a), x\trr b]+(x\trl a)\trr \al_A(b)- (x\trl b)\trr \al_A(a),$
\item[(B4)] $[x, y]\trl \al_A(a)=[\al_V(x), y\trl a]+[x\trl a, \al_V(y)]+\al_V(x)\trl (y\trr a)- \al_V(y)\trl(x\trr a).$
\item[(B5)] $[x, y]\trr \al_A(a)=\al_V(x)\trr (y\trr a)-\al_V(y)\trr (x\trr a)+ \si(\al_V(x), y\trl a)+\si(x\trl a, \al_V(y))+[\al_A(a),\si(x, y)],$
\item[(B6)] $\al_V(x)\trr\si(y,z)+\al_V(y)\trr\si(z,x)+\al_V(z)\trr\si(x,y)=\si([x, y],\al_V(z))+\si([y, z], \al_V(x))+\si([z, x], \al_V(y)),$
\item[(B7)] $[[x, y], \al_V(z)]+[[y, z], \al_V(x)]+[[z, x], \al_V(y)]= \al_V(x)\trl \si(y, z)+\al_V(y)\trl \si(z, x)+\al_V(z)\trl \si(x, y).$
\end{enumerate}
\end{corollary}
In this case, $(A,\al_A)$ is a subalgebra of $(E,\al_E)$ and  $(V,\al_V)$ is in fact a $\sigma$-Hom-Lie algebra acting on $(A,\al_A)$.

Denote the set of all Hom-Lie algebra extending datum of $(A,\al_A)$ by $(V,\al_V)$ of type (a2)  by $\mathcal{A}^{(2)}({A},V)$.

Note that $A{}_{\,\trr}\# {}_{\trl, \sigma}H$ is a Hom-Lie algebra containing $(A,\al_A)$ as a subalgebra.
In fact, any Hom-Lie algebraic structure on $(E,\al_E)$ containing $(A,\al_A)$ as a subalgebra is isomorphic to such a unified product.
\begin{lemma} \label{lem:33-2}
Let $({A},[\cdot,\cdot])$ be a Hom-Lie algebra and $(E,\al_E)$ a Hom-vector space containing $(A,\al_A)$ as a subspace. Suppose that there is a Hom-Lie algebraic structure $(E, [\cdot,\cdot])$ on $(E,\al_E)$ such that $({A},[\cdot,\cdot])$ is a Lie subalgebra of $(E,\al_E)$. Then there exists a Hom-Lie algebraic extending system $\Omega^{(2)}({A},V)$ of $(A,\al_A)$ by $(V,\al_V)$ such that
$(E, [\cdot,\cdot],\delta_E)\cong A{}_{\,\trr}\# {}_{\trl, \sigma}V$.
\end{lemma}

\begin{lemma} \label{l4}
Let $\Omega^{(2)}({A},V)=(\trr, \trl, \sigma, [\cdot,\cdot])$ and ${\Omega}'^{(2)}({A},V)=(\trr', \trl', \sigma', [\cdot,\cdot]')$ be two Hom-Lie algebraic extending datums of $(A,\al_A)$ by $(V,\al_V)$ of type (a2) and $A{}_{\,\trr}\# {}_{\trl, \sigma}V$, $A_{\trr'}\# {}_{\trl', \sigma'}V$ be the corresponding unified products. Then there exists a bijection between the set of all homomorphisms of Hom-Lie algebras $\varphi:A{}_{\,\trr}\# {}_{\trl, \sigma}V\to A_{\trr'}\# {}_{\trl', \sigma'}V$ whose restriction on $(A,\al_A)$ is the identity map and the set of pairs $(r,s)$, where $r:V\rightarrow {A}$ and $s:V\rightarrow V$ are two linear maps satisfying
\begin{eqnarray}
{s}(x)\trl' a&=&{s}(x\trl a),\\
{r}(x\trl a)&=&[{r}(x),a]-x\trr a+ {s}(x)\trr' a,\\
{s}([x, y])&=&[{s}(x),{s}(y)]'+{s}(x)\trl'{r}(y)-{s}(y)\trl'{r}(x),\\
\notag{r}([x, y])&=&[{r}(x),{r}(y)]+{s}(x)\trr'{r}(y)-{s}(y)\trr'{r}(x)+\sigma'({s}(x),{s}(y))-\sigma(x,y)\\
\end{eqnarray}
for all $a\in{A}$ and $x$, $y\in V$.

Under the above bijection the homomorphism of Hom-Lie algebras $\varphi=\varphi_{r,s}: A{}_{\,\trr}\# {}_{\trl, \sigma}V\to A_{\trr'}\# {}_{\trl', \sigma'}V$ to $(r,s)$ is given  by $\varphi(a,x)=(a+r(x), s(x))$ for all $a\in {A}$ and $x\in V$. Moreover, $\varphi=\varphi_{r,s}$ is an isomorphism if and only if $s: V\rightarrow V$ is a linear isomorphism.
\end{lemma}

Let $({A},\al_A)$ be a Hom-Lie algebra and $(V,\al_V)$ a Hom-vector space. Two Hom-Lie algebra extending systems $\Omega^{(i)}({A}, V)$ and ${\Omega'^{(i)}}({A}, V)$  are called equivalent if $\varphi_{r,s}$ is an isomorphism.  We denote it by $\Omega^{(i)}({A}, V)\equiv{\Omega'^{(i)}}({A}, V)$.
From the above lemmas, we obtain the following result.

\begin{theorem}\label{thm3-1}
Let $({A}, [\cdot,\cdot])$ be a Hom-Lie algebra, $(E,\al_E)$ a Hom-vector space containing $(A,\al_A)$ as a subspace and
$(V,\al_V)$ be a complement of $(A,\al_A)$ in $(E,\al_E)$.
Denote $\mathcal{HA}(V,{A}):=\mathcal{A}^{(1)}({A},V)\sqcup \mathcal{A}^{(2)}({A},V) /\equiv$. Then the map
\begin{eqnarray}
&&\Psi: \mathcal{HA}(V,{A})\rightarrow Extd(E,{A}),\\
&&\overline{\Omega^{(1)}({A},V)}\mapsto A_{\theta}\#_{\trl} V,\quad \overline{\Omega^{(2)}({A},V)}\mapsto A{}_{\,\trr}\# {}_{\trl, \sigma} V
\end{eqnarray}
is bijective, where $\overline{\Omega^{(i)}({A}, V)}$ is the equivalence class of $\Omega^{(i)}({A}, V)$ under $\equiv$.
\end{theorem}

\subsection{Extending structures for Hom-Lie coalgebras}
Next we consider the Hom-Lie coalgebra structures on $E=A^{\phi, P}\# {}^{\psi, Q} V$.

There are two cases for $(A,\delta_A)$ to be a Hom-Lie coalgebra. The first case is when $\phi\neq 0, Q=0$,  we obtain the following type (c1) unified product for Hom-Lie coalgebras.
\begin{corollary}\label{cor03}
Let $({A},\delta_{A})$ be a Hom-Lie coalgebra and $(V,\al_V)$ a Hom-vector space.
An  extending datum  of $(A,\al_A)$ by $(V,\al_V)$ of  type (c1) is  $\Omega^c({A},V)=(\phi, {\psi}, P, \delta_V)$ with  linear maps
\begin{eqnarray*}
 \phi :A \to V\otimes A, \quad{\psi}: V\to  V\otimes A,~~~~{P}: A\rightarrow {V}\otimes {V},~~~~\delta_V: V\rightarrow V\otimes V.
\end{eqnarray*}
 Denote by $A^{\phi, P}\# {}^{\psi} V$ the vector space $E={A}\oplus V$ with the linear map
$\delta_E: E\rightarrow E\otimes E$ given by
\begin{eqnarray}\label{eq11}
\delta_E(a)=\delta_{A}(a)+{\phi}(a)-\tau{\phi}(a)+P(a),\quad
\delta_E(x)=\delta_V(x)+{\psi}(x)-\tau{\psi}(x).
\end{eqnarray}
Then $A^{\phi, P}\# {}^{\psi} V$  is a Hom-Lie coalgebra with the Lie cobracket given by (\ref{eq11}) if and only if the following compatibility conditions hold:
\begin{enumerate}
\item[(C1)] $ P(a)=-\tau P(a),~~~\delta_V(x)=-\tau\delta_V(x),$
\item[(C2)] $\delta_V(a\lmoi)\ot \al_A(a\lmoo)+P(a\li)\ot  \al_A(a\lii)=\al_V(a\lmoi)\ot \phi(a\lmoo) - \tau_{12}\left(a\lmoi\ot \phi(a\lmoo)\right)$\\
$+\al_V(a\ppi)\ot \psi\left(a\pii\right)  +\tau_{23}\left(\psi(a\ppi)\ot \al_V(a\pii)\right),$
\item[(C3)] $\al_V(x\lmoo)\ot \delta_A(x\lmi) = \psi(x\lmoo)\ot \al_A(x\lmi) -  \tau_{23}\left(\psi(x\lmoo)\ot \al_A(x\lmi)\right),$
\item[(C4)] $\al_V(a\loi)\ot \delta_A(a\loo) =  \phi(a\li)\ot \al_A(a\lii)+  \tau_{12}\left( \al_A(a\li)\ot \phi(a\lii)\right)$\\
$+\psi(a\loi)\ot \al_A(a\loo)-  \tau_{23}\left( \psi(a\loi)\ot \al_A(a\loo)\right),$
\item[(C5)] $ \delta_V(x\lmo)\ot \al_A(x\lmi)= \al_V(x\li)\ot \psi(x\lii) +  \tau_{23}\left( \psi(x\li)\ot \al_V(x\lii)\right)$\\
$+\al_V(x\lmo)\ot \psi(x\lmi) -  \tau_{12}\left( \al_V(x\lmo)\ot \psi(x\lmi)\right),$
\item[(C6)] $\al_V(a\lmoi)\ot P(a\lmoo)+\tau_{12}\tau_{23}\left(\al_V(a\lmoi)\ot P(a\lmoo)\right)+\tau_{23}\tau_{12}\left(\al_V(a\lmoi)\ot P(a\lmoo)\right)$\\
$=\delta(a\ppi)\ot \al_V(a\pii)+\tau_{12}\tau_{23}\left(\delta(a\ppi)\ot \al_V(a\pii)\right)+\tau_{23}\tau_{12}\left(\delta(a\ppi)\ot \al_V(a\pii)\right),$
\item[(C7)] $\delta(x\li)\ot \al_V(x\lii)+\tau_{12}\tau_{23}\left(\delta(x\li)\ot \al_V(x\lii)\right)+\tau_{23}\tau_{12}\left( \delta(x\li)\ot \al_V(x\lii)\right)$\\
$= \al_V(x\lmoo)\ot P(x\lmi)+\tau_{12}\tau_{23}\left(\al_V(x\lmoo)\ot P(x\lmi)  \right)+\tau_{23}\tau_{12}\left(\al_V(x\lmoo)\ot P(x\lmi) \right).$
 \end{enumerate}
\end{corollary}
Denote the set of all Hom-Lie coalgebra extending datum of $(A,\al_A)$ by $(V,\al_V)$ of type (c1) by $\mathcal{C}^{(1)}({A},V)$.

In this case, although $(A,\al_A)$ is a Hom-Lie coalgebra but it is not a subcoalgebra of  $A^{\phi, P}\# {}^{\psi} V$.
The characterization of this type of Hom-Lie coalgebra $(A,\al_A)$ is as follows.
\begin{lemma}\label{lem:33-3}
Let $({A},\delta_{A})$ be a Hom-Lie coalgebra and $(E,\al_E)$ a Hom-vector space containing $(A,\al_A)$ as a subspace. Suppose that there is a Hom-Lie coalgebra structure $(E,\delta_E)$ on $(E,\al_E)$ such that  $p: E\to {A}$ is a Hom-Lie coalgebra homomorphism. Then there exists a Hom-Lie coalgebra extending system $\Omega^c({A}, V)$ of $({A},\delta_{A})$ by $(V,\al_V)$ such that $(E,\delta_E)\cong A^{\phi, P}\# {}^{\psi} V$.
\end{lemma}

\begin{proof}
Let $p: E\to {A}$ and $\pi: E\to V$ be the projection maps and $V=\ker({p})$.
Then the extending datum of $({A},\delta_{A})$ by $(V,\al_V)$ is defined as follows:
\begin{eqnarray*}
&&{\phi}: A\rightarrow V\ot {A},~~~~{\phi}(x)=(\pi\otimes {p})\delta_E(a),\\
&&{\psi}: V\rightarrow V\ot {A},~~~~{\phi}(x)=(\pi\otimes {p})\delta_E(x),\\
&&\delta_V: V\rightarrow V\otimes V,~~~~\delta_V(x)=(\pi\otimes \pi)\delta_E(x),\\
&&P: A\rightarrow {V}\otimes {V},~~~~P(a)=({\pi}\otimes {\pi})\delta_E(a).
\end{eqnarray*}
One check that  $\varphi: A^{\phi, P}\# {}^{\psi} V\to E$ given by $\varphi(a,x)=a+x$ for all $a\in A, x\in V$ is a Hom-Lie coalgebra isomorphism.
\end{proof}

\begin{lemma}\label{lem-c1}
Let $\Omega^{(1)}({A}, V)=(\phi, {\psi}, P, \delta_V)$ and ${\Omega'^{(1)}}({A}, V)=(\phi', {\psi'}, P', \delta'_V)$ be two Hom-Lie coalgebra extending datums of $({A},\delta_{A})$ by $(V,\al_V)$. Then there exists a bijection between the set of  Hom-Lie coalgebra homomorphisms $\varphi: A^{\phi, P}\# {}^{\psi} V\rightarrow A^{\phi', P'}\# {}^{\psi'} V$ whose restriction on $(A,\al_A)$ is the identity map and the set of pairs $(r,s)$, where $r:V\rightarrow {A}$ and $s:V\rightarrow V$ are two linear maps satisfying
\begin{eqnarray}
\label{comorph11}&&P'(a)=s(a\ppi)\ot s(a\pii),\\
\label{comorph12}&&\phi'(a)={s}(a\lmoi)\ot a\lmo+s(a\ppi)\ot r(a\pii),\\
\label{comorph13}&&\delta'_A(a)=\delta_A(a)+{r}(a\lmoi)\ot a\lmo-a\lmo\ot {r}(a\lmoi)+r(a\ppi)\ot r(a\pii)\\
\label{comorph21}&&\delta_V'({s}(x))=({s}\otimes {s})\delta_V(x),\\
\label{comorph22}&&{\psi}'({s}(x))=s(x\li)\ot r(x\lii)+s(x\lmo)\ot x\lmi,\\
\label{comorph23}&&\delta'_A({r}(x))=r(x\li)\ot r(x\lii)+r(x\lmo)\ot x\lmi-x\lmi\ot r(x\lmo).
\end{eqnarray}
Under the above bijection the Hom-Lie coalgebra homomorphism $\varphi=\varphi_{r,s}: A^{\phi, P}\# {}^{\psi} V\rightarrow A^{\phi', P'}\# {}^{\psi'} V$ to $(r,s)$ is given by $\varphi(a,x)=(a+r(x),s(x))$ for all $a\in {A}$ and $x\in V$. Moreover, $\varphi=\varphi_{r,s}$ is an isomorphism if and only if $s: V\rightarrow V$ is a linear isomorphism.
\end{lemma}

In the case $\phi=0, Q\neq 0$, then from (TBB3) we get that $a\ppi\ot Q(a\pii)=0$, since $Q\neq 0$ we assume $P=0$ for simplicity, thus we obtain  the following type (c2) unified product for Hom-Lie coalgebras.
\begin{corollary}\label{cor:04}
Let $({A},\delta_{A})$ be a Hom-Lie coalgebra and $(V,\al_V)$ a Hom-vector space.
An  extending datum  of $({A},\delta_{A})$ by $(V,\al_V)$ of type (c2)  is  $\Omega^{(2)}({A},V)=({\psi}, {Q}, \delta_V)$ with  linear maps
\begin{eqnarray*}
{\psi}: V\to  V\otimes A,~~~~{Q}: V\rightarrow {A}\otimes {A},~~~~\delta_V: V\rightarrow V\otimes V.
\end{eqnarray*}
 Denote by $A^{}\# {}^{\psi, Q} V$ the vector space $E={A}\oplus V$ with the linear map
$\delta_E: E\rightarrow E\otimes E$ given by
\begin{eqnarray}\label{eq1}
\delta_E(a)=\delta_{A}(a),\quad
\delta_E(x)=\delta_V(x)+{\psi}(x)-\tau{\psi}(x)+Q(x).
\end{eqnarray}
Then $A^{}\# {}^{\psi, Q} V$  is a Hom-Lie coalgebra with the Lie cobracket given by (\ref{eq1}) if and only if the following compatibility conditions hold:
\begin{enumerate}
\item[(D1)] $Q(x)=-\tau Q(x),~~~\delta_V(x)=-\tau\delta_V(x),$
\item[(D2)] $\delta_V(x\li)\ot \al_V(x\lii)+\tau_{12}\tau_{23}\left(\delta_V(x\li)\ot \al_V(x\lii)\right)+\tau_{23}\tau_{12}\left( \delta_V(x\li)\ot \al_V(x\lii) \right)=0,$
\item[(D3)] $\al_V(x\lmoo)\ot \delta_A(x\lmi)+\al_V(x\li)\ot \psi(x\lii)= \psi(x\lmoo)\ot \al_A(x\lmi) -  \tau_{23}\left(\psi(x\lmoo)\ot \al_A(x\lmi)\right)$\\
   $+\psi(x\qi)\ot \al_A(x\qii) + \tau_{12}\left(\al_A(x\qi)\ot\psi( x\qii)\right),$
\item[(D4)] $\delta_V(x\lmo)\ot \al_A(x\lmi) =  \al_V(x\li)\ot \psi(x\lii) +  \tau_{23}\left( \psi(x\li)\ot \al_V(x\lii)\right),$
\item[(D5)] $Q(x\lmoo)\ot \al_A(x\lmi)+\tau_{12}\tau_{23}\left(Q(x\lmoo)\ot \al_A(x\lmi)\right)+\tau_{23}\tau_{12}\left(Q(x\lmoo)\ot \al_A(x\lmi)\right)$\\
 $= \al_A(x\qi)\ot \delta_A(x\qii)+\tau_{12}\tau_{23}\left(\al_A(x\qi)\ot \delta_A(x\qii) \right)+\tau_{23}\tau_{12}\left(\al_A(x\qi)\ot \delta_A(x\qii)\right).$
\end{enumerate}
\end{corollary}
Denote the set of all Hom-Lie coalgebra extending datum of $(A,\al_A)$ by $(V,\al_V)$ of type (c2) by $\mathcal{C}^{(2)}({A},V)$.

Similar as Hom-Lie algebra case,  one  show that any Hom-Lie coalgebra structure on $(E,\al_E)$ containing $(A,\al_A)$ as a subcoalgebra is isomorphic to such a unified coproduct.
\begin{lemma}\label{lem:33-4}
Let $({A},\delta_{A})$ be a Hom-Lie coalgebra and $(E,\al_E)$ a Hom-vector space containing $(A,\al_A)$ as a subspace. Suppose that there is a Hom-Lie coalgebra structure $(E,\delta_E)$ on $(E,\al_E)$ such that  $({A},\delta_{A})$ is a Lie subcoalgebra of $(E,\al_E)$. Then there exists a Hom-Lie coalgebra extending system $\Omega^{(2)}({A}, V)$ of $({A},\delta_{A})$ by $(V,\al_V)$ such that $(E,\delta_E)\cong A^{}\# {}^{\psi, Q} V$.
\end{lemma}

\begin{proof}
Let $p: E\to {A}$ and $\pi: E\to V$ be the projection map and $V=ker({p})$.
Then the extending datum of $({A},\delta_{A})$ by $(V,\al_V)$ is defined as follows:
\begin{eqnarray*}
&&{\psi}: V\rightarrow V\ot {A},~~~~{\phi}(x)=(\pi\otimes {p})\delta_E(x),\\
&&\delta_V: V\rightarrow V\otimes V,~~~~\delta_V(x)=(\pi\otimes \pi)\delta_E(x),\\
&&Q: V\rightarrow {A}\otimes {A},~~~~Q(x)=({p}\otimes {p})\delta_E(x).
\end{eqnarray*}
One check that  $\varphi: A^{}\# {}^{\psi, Q} V\to E$ given by $\varphi(a,x)=a+x$ for all $a\in A, x\in V$ is a Hom-Lie coalgebra isomorphism.
\end{proof}

\begin{lemma}\label{lem-c2}
Let $\Omega^{(2)}({A}, V)=( {\psi}, {Q}, \delta_V)$ and ${\Omega'^{(2)}}({A}, V)=({\psi'}, {Q'}, \delta'_V)$ be two Hom-Lie coalgebra extending datums of $({A},\delta_{A})$ by $(V,\al_V)$. Then there exists a bijection between the set of  Hom-Lie coalgebra homomorphisms $\varphi: A \# {}^{\psi, Q} V\rightarrow A \# {}^{\psi', Q'} V$ whose restriction on $(A,\al_A)$ is the identity map and the set of pairs $(r,s)$, where $r:V\rightarrow {A}$ and $s:V\rightarrow V$ are two linear maps satisfying
\begin{eqnarray}
\label{comorph2}&&{\psi}'({s}(x))=s(x\li)\ot r(x\lii)+s(x\lmo)\ot x\lmi,\\
\label{comorph3}&&\delta_V'({s}(x))=({s}\otimes {s})\delta_V(x),\\
\label{comorph4}&&\delta'_A({r}(x))+{Q'}({s}(x))=r(x\li)\ot r(x\lii)+r(x\lmo)\ot x\lmi-x\lmi\ot r(x\lmo)+{Q}(x).
\end{eqnarray}
Under the above bijection the Hom-Lie coalgebra homomorphism $\varphi=\varphi_{r,s}: A^{ }\# {}^{\psi, Q} V\rightarrow A^{ }\# {}^{\psi', Q'} V$ to $(r,s)$ is given by $\varphi(a,x)=(a+r(x),s(x))$ for all $a\in {A}$ and $x\in V$. Moreover, $\varphi=\varphi_{r,s}$ is an isomorphism if and only if $s: V\rightarrow V$ is a linear isomorphism.
\end{lemma}

Let $({A},\delta_{A})$ be a Hom-Lie coalgebra and $(V,\al_V)$ a Hom-vector space. Two Hom-Lie coalgebra extending systems $\Omega^{(i)}({A}, V)$ and ${\Omega'^{(i)}}({A}, V)$  are called equivalent if $\varphi_{r,s}$ is an isomorphism.  We denote it by $\Omega^{(i)}({A}, V)\equiv{\Omega'^{(i)}}({A}, V)$.
From the above lemmas, we obtain the following result.
\begin{theorem}\label{thm3-2}
Let $({A},\delta_{A})$ be a Hom-Lie coalgebra, $(E,\al_E)$ a Hom-vector space containing $(A,\al_A)$ as a subspace and
$(V,\al_V)$ be a ${A}$-complement in $(E,\al_E)$. Denote $\mathcal{HC}(V,{A}):=\mathcal{C}^{(1)}({A},V)\sqcup\mathcal{C}^{(2)}({A},V) /\equiv$. Then the map
\begin{eqnarray}
&&\Psi: \mathcal{HC}_{{A}}^2(V,{A})\rightarrow CExtd(E,{A}),\\
&&\overline{\Omega^{(1)}({A},V)}\mapsto A^{\phi, P}\# {}^{\psi} V, \quad \overline{\Omega^{(2)}({A},V)}\mapsto A^{}\# {}^{\psi, Q} V
\end{eqnarray}
is bijective, where $\overline{\Omega^{(i)}({A},V)}$ is the equivalence class of $\Omega^{(i)}({A}, V)$ under $\equiv$.
\end{theorem}


%

\subsection{Extending structures  for Hom-Lie bialgebras}

There are two special cases for which $({A}, [\cdot,\cdot], \delta_{A})$ is reduced to a Hom-Lie bialgebra.
The first case is when $\trr=0, \sigma=0, Q=0$ in the above Theorem \ref{main2}. In this case we obtain the following result.

\begin{theorem}\label{thm-41}
Let $(A,[\cdot,\cdot],\delta_A)$ be a Hom-Lie bialgebra and $(V,\al_V)$ a Hom-vector space.
An extending datum of $(A,\al_A)$ by $(V,\al_V)$ of type (I) is  $\Omega^{(1)}({A},V)=(\trl,  \phi, \psi, P, [\cdot,\cdot]_V, \delta_V)$ consisting of  linear maps
\begin{eqnarray*}
\trl: V\times {A}\rightarrow V,~~~~\theta:  A\times A \rightarrow {V},~~~[\cdot,\cdot]_V:V\times V \rightarrow V,\\
 \phi :A \to V\otimes A, \quad{\psi}: V\to  V\otimes A,~~~~{P}: A\rightarrow {V}\otimes {V},~~~~\delta_V: V\rightarrow V\otimes V.
\end{eqnarray*}
Then the unified product $A^{\phi,P}_{\theta}\# {}^{\psi}_{\trl}\, V$ with bracket
\begin{align}
[(a, x), (b, y)]:=([a, b], [x,y]+x\trl b-y\trl a+\theta(a, b))
\end{align}
and cobracket
\begin{eqnarray}
\delta_E(a)=\delta_{A}(a)+{\phi}(a)-\tau{\phi}(a)+P(a),\quad \delta_E(x)=\delta_V(x)+{\psi}(x)-\tau{\psi}(x)
\end{eqnarray}
form a Hom-Lie bialgebra if and only if $A_{\theta}\# {}_{\trl} V$ form a Hom-Lie algebra, $A^{\phi,P}\# {}^{\psi} \, V$ form a Hom-Lie coalgebra and the following conditions are satisfied:
\begin{enumerate}
\item[(E1)] $ \delta_{V} (x\trl a)=\al_V(x\li)\ot x\lii \trl \al_A(a)+x\li\trl \al_A(a)\ot x\lii$\\
$+ \al_V(a\loi)\ot \al_V(x)\trl a\loo-\al_V(x)\trl a\loo\ot \al_V(a\loi)$\\
$+[\al_V(x), P(a)]+\al_V(x\lmoo)\ot \theta(x\lmi, \al_A(a))-\theta(x\lmi, \al_A(a))\ot \al_V(x\lmoo)$,
\item[(E2)] $\phi([a, b])+\psi\theta(a, b)= \al_V(a\lmoi)\ot [a\lmoo, \al_A(b)]+\al_H(b\lmoi)\ot [\al_A(a), b\lmoo]$\\
$+a\lmoi\trl \al_A(b) \ot \al_A(a\lmoo)-b\lmoi\trl \al_A(a) \ot \al_A(b\lmoo)$\\
$+\theta(\al_A(a), b\li)\ot \al_A(b\lii)+\theta(a\li, \al_A(b))\ot \al_A(a\lii)$,

\item[(E3)] $\psi([x, y])= [\al_V(x), y\lmoo]\ot \al_A(y\lmi)+[x\lmoo, \al_V(y)]\ot \al_A(x\lmi)$,
\item[(E4)]
$\delta_V\theta(a, b)+P([a,b])= \al_V(a\lmoi)\ot\theta(a\lmoo, \al_A(b))+\al_H(b\lmoi)\ot\theta(\al_A(a),b\lmoo)$\\
$- \theta(\al_A(a), b\lmoo)\ot \al_H(b\lmoi) -\theta(a\lmoo, \al_A(b))\ot \al_V(a\lmoi)$\\
$+\al_V(a\ppi)\ot a\pii \trl \al_A(b)+a\ppi\trl \al_A(b) \ot \al_V(a\pii)$\\
$-b\ppi\trl \al_A(a)\ot \al_V(b\pii) -\al_V(b\ppi)\ot b\pii\trl \al_A(a)$,
\item[(E5)] $\delta_V([x, y])= [\delta_V(x), \al_V(y)]+[\al_V(x), \delta_V([y)]$\\
$-x\lmoo\ot \al_V(y)\trl x\lmi -\al_V(x)\trl y\lmi\ot y\lmoo$\\
$+ y\lmoo\ot \al_V(x)\trl y\lmi+\al_V(y)\trl x\lmi\ot x\lmoo$,
\item[(E6)] $ \psi (x\trl a)=\al_V(x\loo)\ot [x\lmi, \al_A(a)]+x\loo\trl \al_A(a)\ot \al_A(x\lmi) +\al_V(x)\trl a\li\ot \al_A(a\lii),$
\end{enumerate}
Conversely, any Hom-Lie bialgebra structure on $(E,\al_E)$ with the canonical projection map $p: E\to A$ both a Hom-Lie algebra homomorphism and a Hom-Lie coalgebra homomorphism is of this form.
\end{theorem}
Note that in this case, although $(A,[\cdot,\cdot],\delta_A)$ is not a Lie sub-bialgebra of $A^{\phi,P}_{}\# {}^{\psi,Q}_{\trl, \sigma}\, V$, but it is indeed a Hom-Lie bialgebra and a subspace $A^{\phi,P}_{}\# {}^{\psi,Q}_{\trl, \sigma}\, V$.
Denote the set of all Hom-Lie bialgebra extending datum of type (I) by $\mathcal{LB}^{(1)}({A},V)$.

The second case is when $\theta=0,  P=0, \phi=0$ in the above Theorem \ref{main2}. In this case we obtain the following result.
\begin{theorem}\label{thm-42}
Let $(A,\al_A)$ be a Hom-Lie bialgebra and $(V,\al_V)$ a Hom-vector space.
An extending datum of $(A,\al_A)$ by $(V,\al_V)$ of type (II) is  $\Omega^{(2)}({A},V)=(\trr, \trl, \sigma, \psi, Q,  [\cdot,\cdot]_V, \delta_V)$ consisting of  linear maps
\begin{eqnarray*}
\trr: V\times {A}\rightarrow {A},~~~~\trl: V\times {A}\rightarrow V,~~~~\sigma:  V\times V \rightarrow {A},~~~[\cdot,\cdot]_V:V\times V \rightarrow V,\\
{\psi}: V\to  V\otimes A,~~~~{Q}: V\rightarrow {A}\otimes {A},~~~~\delta_V: V\rightarrow V\otimes V.
\end{eqnarray*}
Then the unified product $A^{}_{\trr}\# {}^{\psi,Q}_{\trl, \sigma}\, V$ with bracket
\begin{align}
[(a, x), (b, y)]:=([a, b]+x\trr b-y\trr a+\sigma(x, y), [x,y]+x\trl b-y\trl a)
\end{align}
and cobracket
\begin{eqnarray}
\delta_E(a)=\delta_{A}(a),\quad \delta_E(x)=\delta_V(x)+{\psi}(x)-\tau{\psi}(x)+Q(x)
\end{eqnarray}
form a Hom-Lie bialgebra if and only if $A{}_{\,\trr}\# {}_{\trl, \sigma} V$ form a Hom-Lie algebra, $A^{}\# {}^{\psi,Q} \, V$ form a Hom-Lie coalgebra and the following conditions are satisfied:
\begin{enumerate}
\item[(F1)] $\delta_{A} (x\trr a)+Q(x\trl a)=\al_V(x)\trr a\li\ot \al_A(a\lii)+\al_A(a\li)\ot \al_V(x)\trr a\lii$\\
$\qquad+ x\loo\trr \al_A(a)\ot \al_A(x\lmi)-\al_A(x\lmi)\ot x\loo\trr\al_A(a)+[Q(x), \al_A(a)]$,
\item[(F2)] $ \delta_V(x\trl a)=\al_V(x\li)\ot x\lii \trl \al_A(a)+x\li\trl \al_A(a)\ot x\lii$,
\item[(F3)] $  \psi([x, y])=[\al_V(x), y\loo]\ot \al_A(y\lmi)+[x\loo, \al_V(y)]\ot \al_A(x\lmi)$\\
$+\al(y\loo)\ot \al_V(x)\trr y\lmi -\al_V(x\loo)\ot \al_V(y)\trr x\lmi$\\
$+ \al_V(x)\li\ot\si(x\lii, \al_V(y))+\al_V(y\li)\ot \si(\al_V(x), y\lii)$\\
$+\al_V(x)\trl y\qi\ot \al_A(y\qii)-\al_V(y)\trl x\qi\ot \al_A(x\qii)$,
\item[(F4)] $\delta_A \si(x, y)+Q([x, y])=\si(x\lmoo,  y)\ot \al_A(x\lmi)+\si(\al_V(x), y\lmoo)\ot \al_A(y\lmi)$\\
$-\al_A(y\lmi)\ot \si(\al_V(x), y\lmoo)-\al_A(x\lmi)\ot \si(x\lmoo, y)$\\
$+\al_V(x)\trr y\qi\ot \al_A(y\qii)+\al_A(y\qi)\ot \al_V(x)\trr y\qii $\\
$-\al_A(x\qi)\ot  \al_V(y)\trr  x\qii - \al_V(y)\trr x\qi \ot \al_A(x\qii)$,
\item[(F5)] $\delta_V([x, y]) = [\delta_V(x), \al_V(y)]+[\al_V(x), \delta_V([y)] $\\
$-y_{(-1)}  \trr \alpha (x) \otimes \al_V(y_{(0)})-\al_V(x_{(0)}) \otimes  x_{(-1)}  \trr \alpha_V(y)$\\
$+x_{(-1)} \trr \alpha_V(y) \otimes \al_V(x_{(0)})+\al_V(y_{(0)}) \otimes  y_{(-1)}  \trr \alpha_V(x)$,
\item[(F6)] $\psi (x\trl a)=\al_V(x\li) \otimes x_{2} \trr \alpha_A(a)$\\
 $+\al_V(x\loo)\ot [x\lmi, \al_A(a)]+x\loo\trl \al_A(a)\ot \al_A(x\lmi) +\al_V(x)\trl a\li\ot \al_A(a\lii)$.
\end{enumerate}
Conversely, any Hom-Lie bialgebra structure on $(E,\al_E)$ with the canonical injection map $i: A\to E$ both a Hom-Lie algebra homomorphism and a Hom-Lie coalgebra homomorphism is of this form.
\end{theorem}
Denote the set of all  Hom-Lie bialgebra extending datum of type (II) by $\mathcal{LB}^{(2)}({A},V)$.

Note that $A^{\phi,P}_{}\# {}^{\psi}_{\trl, \sigma}\, V$ and $A^{}_{\trr}\# {}^{\psi,Q}_{\trl, \sigma}\, V$ are all Hom-Lie bialgebra structures on $(E,\al_E)$.
Conversely,  any Hom-Lie bialgebra extending system $(E,\al_E)$ of $(A,\al_A)$  through $(V,\al_V)$ is isomorphic to such a unified products of the two types.
Now from Theorem \ref{thm3-1}, Theorem \ref{thm3-2}  in last section and Theorem \ref{thm-41}, Theorem \ref{thm-42} we obtain the main result of in this section,
which solve the extending problem for Hom-Lie bialgebra.


\begin{theorem}\label{thm43}
Let $({A}, [\cdot,\cdot], \delta_{A})$ be a Hom-Lie bialgebra, $(E,\al_E)$ a Hom-vector space containing $(A,\al_A)$ as a subspace and $(V,\al_V)$ be a complement of $(A,\al_A)$ in $(E,\al_E)$.
Denote by
$$\mathcal{HLB}(V,{A}):=\mathcal{LB}^{(1)}({A},V)\sqcup\mathcal{LB}^{(2)}({A},V)/\equiv.$$
Then the map
\begin{eqnarray}
&&\Upsilon: \mathcal{HLB}(V,{A})\rightarrow BExtd(E,{A}),\\
&&\overline{\Omega^{(1)}({A},V)}\mapsto A^{\phi,P}_{\theta}\# {}^{\psi}_{\trl}\, V,\quad   \overline{\Omega^{(2)}({A},V)}\mapsto A^{}_{\trr}\# {}^{\psi,Q}_{\trl, \sigma}\, V
\end{eqnarray}
is bijective, where $\overline{\Omega^{(i)}({A}, V)}$ is the equivalence class of $\Omega^{(i)}({A}, V)$ under $\equiv$.
\end{theorem}

\subsection{Flag extending structures}
In this section, we study the case when $(V,\al_V)$ is a 1-dimensional vector space. This will be  called flag extending system.
Since $(V,\al_V)$ is a 1-dimensional vector space, then the bracket and cobracket of $(V,\al_V)$ is given by $[x,y]=0$ and $\delta_V(x)=0$ for all $x,y\in V$.

\begin{lemma}
Let $({A},[\cdot,\cdot],\delta_{{A}})$ be a braided Hom-Lie bialgebra and $V=k\{x\}$ be a 1-dimensional vector space. A  flag datum consists of
$$\lambda: A\to k,\quad D:A\to A, \quad T:A\to A, \quad a_0\in A$$
satisfying the following compatibility conditions:
\begin{eqnarray}
&&\lambda([a,b])=\lam(\al_A(a))\lam(b)-\lam(\al_A(b))\lam(a),\\
&&D([a,b])=[D(a),\al_A(b)]+[\al_A(a),D(b)]+\lambda(a)D(b)-\lambda(b)D(a),\\
&&T([a, b])=[T(a),\al_A(b)]+[\al_A(a),T(b)]+\lambda(b)T(a)-\lambda(a)T(b),\\
&&T(D(a))=D(T(a))+[a_0,\al_A(a)]+\lambda(a_1)\al(a_2).
\end{eqnarray}
The corresponding the extending datum $\Omega(A, V)$  is given by:
\begin{eqnarray}
&&x\triangleright a =D(a), \quad x\triangleleft a=  \lambda(a)x,  \quad \phi(a)= x\otimes T(a),\quad \psi(x)=x\otimes a_0,\\
&&\sigma(x, x) = 0, \quad  [x, x] = 0,\quad P(a)=0,\quad \delta_V(x)=0.
\end{eqnarray}
The unified product associated to this flag extending system is given by
\begin{equation}
[(a, x), (b, y)] =\Big ([a, b] +  D(a)y -  D(b)x , \lambda (a)y - \lambda (b)) x\Big),
\end{equation}
and
\begin{eqnarray}
&&\delta_E(a)=\delta_A(a)+ x\otimes T(a)-T(a)\otimes x,\quad \delta_E(x)=x\otimes a_0-a_0\otimes x.
\end{eqnarray}
\end{lemma}

Denote the set of all flag datums of braided Hom-Lie bialgebra by $\mathcal{FB}({A})$.

\begin{Definition}
Two flag datums $(\lambda, D, T, a_0)$ and $(\lambda', D', T, a_0') \in \mathcal{FB}({A})$ are called equivalent if $\lambda'=\lambda$, $a_0'=a_0$ and there exist some element $r_0\in {A}$  such that
\begin{eqnarray}
\label{bbb11}&&\lambda(a)r_0=[r_0,a] - D(a)+s D'(a),\\
\label{bbb12}&&\delta'_{{A}}(r_0)=r_0\ot a_0-a_0\ot r_0,\\
\label{bbb13}&&\delta'_{{A}}(a)=\delta_{{A}}(a)+r_0\otimes T(a)-T(a)\otimes r_0.
\end{eqnarray}
\end{Definition}

By the above lemma, we have
\begin{theorem}
Let $({A},[\cdot, \cdot], \delta_{{A}})$ be a  braided Hom-Lie bialgebra and $(V,\al_V)$ be a 1-dimensional vector space. Then there is a bijection between the set $\mathcal{BLB}({A},V)$ of all  Hom-Lie bialgebra extending systems of $(A,\al_A)$ by $(V,\al_V)$  and $\mathcal{FB}({A})$.
\end{theorem}


Next, we consider flag extending systems for Hom-Lie bialgebras.
\begin{lemma}
Let $({A},[\cdot,\cdot],\delta_{{A}})$ be a Hom-Lie bialgebra. A  flag datum of type (I) consists of
$$\lambda: A\to k,\quad T:A\to A, \quad a_0\in A$$
satisfying the following compatibility conditions:
\begin{eqnarray}
\label{f11}
&&\lambda([a,b])=\lam(\al_A(a))\lam(b)-\lam(\al_A(b))\lam(a),\\
\label{f12}&&T([a, b])=[T(a),\al_A(b)]+[\al_A(a),T(b)],\\
&&[a_0,\al_A(a)]+\lambda(a_1)\al(a_2)=0.\label{f13}
\end{eqnarray}
The corresponding the extending datum $\Omega^{(1)}(A, V)$ of type (I)   is given by:
\begin{eqnarray}
&&x\triangleleft a=  \lambda(a)x,  \quad \phi(a)= x\otimes T(a),\quad \psi(x)=x\otimes a_0,\\
&&\theta(a, b) = 0, \quad  [x, x] = 0,\quad P(a)=0,\quad \delta_V(x)=0.
\end{eqnarray}
The unified product $A\#^{(1)} V$ associated to the flag extending system is given by
\begin{equation}
[(a, x), (b, y)] =\Big ([a, b], \lambda (b)x - \lambda (a)) y\Big),
\end{equation}
and
\begin{eqnarray}
&&\delta_E(a)=\delta_A(a)+ x\otimes T(a)-T(a)\otimes x,\quad \delta_E(x)=x\otimes a_0-a_0\otimes x.
\end{eqnarray}
\end{lemma}

Denote the set of all flag datums of type (I)  by $\mathcal{F}^{(1)}({A})$.

\begin{lemma}
Let $({A},[\cdot,\cdot], \delta_{{A}})$ be a Hom-Lie bialgebra. A  flag datum of type (II) consists of
$$\lambda: A\to k,\quad D:A\to A,\quad a_0\in A,\quad Q\in A\wedge A$$
 satisfying the following compatibility conditions:
\begin{eqnarray}
\label{td1}&&\lambda([a,b])=\lam(\al_A(a))\lam(b)-\lam(\al_A(b))\lam(a),\\
&&\label{td2}D([a,b])=[D(a),\al_A(b)]+[\al_A(a),D(b)]+\lambda(a)D(b)-\lambda(b)D(a),\\
\label{tdctd1}&&[a_0,\al_A(a)]+\lambda(a_1)\al(a_2)=0,\\
\notag&&\delta_{{A}}(D(a))+\lambda(a)Q=[\al_A(a),Q]+qD(a_{1})\otimes \al(a_{2})+q\al(a_{1})\otimes D(a_{2})\\
\label{tdctd2}&&+D(\al_A(a))\otimes \al_A(a_0)-\al_A(a_0)\otimes D(\al_A(a)),\\
\label{tdctd3}&&\al(a_{0})\otimes Q-\tau_{12}\left(\al(a_{0}) \otimes Q\right)+Q\otimes\al(a_{0})=\left(\al \otimes \delta-\tau_{12}(\al \otimes \delta)-\delta \otimes \al\right) Q.
\end{eqnarray}
The corresponding is the extending datum $\Omega^{(2)}(A, V)$ of type (II)   given by:
\begin{eqnarray}
&&x\triangleleft a=  \lambda (a) x,   \quad x\triangleright a =D(a),\quad \omega(x, x) = 0, \\
&& {\psi}(x)=x\otimes a_0,~~~Q(x)=Q.
\end{eqnarray}
The unified product $A\#^{(2)}  V$ is given
\begin{equation}
[(a, x), (b, y)] =\Big ([a, b] +  D(a)y -  D(b)x , \lambda (a)y - \lambda (b)) x\Big).
\end{equation}
and
\begin{equation}
\delta_E(a)=\delta_A(a), \quad \delta_E(x)= x\otimes a_0-a_0 \otimes x+Q.
\end{equation}
\end{lemma}
Denote the set of all flag datums  of type (II)  by $\mathcal{F}^{(2)}({A})$.

By the above two lemmas, we have
\begin{theorem}
Let $({A},[\cdot, \cdot], \delta_{{A}})$ be a Hom-Lie bialgebra and $V=k\{x\}$ be a 1-dimensional vector space. Then there is a bijection between the set $\mathcal{LB}({A},V)$ of all Hom-Lie bialgebra extending systems of $(A,\al_A)$ by $(V,\al_V)$  and $\mathcal{F}({A})=\mathcal{F}^{(1)}({A})\sqcup \mathcal{F}^{(2)}({A})$.
\end{theorem}

\begin{Definition}
Two flag datums $(\lambda, T, a_0)$ and $(\lambda', T', a_0') \in \mathcal{F}^{(1)}({A})$ are called equivalent if $\lambda'=\lambda$, $a_0'=a_0$ and there exist some element $r_0=r(x)\in {A}$ such that
\begin{eqnarray}
\label{eqq11}&&\lambda(a)r_0=[r_0,a],\\
\label{eqq12}&&\delta'_{{A}}(r_0)=r_0\ot a_0-a_0\ot r_0,\\
\label{eqq13}&&\delta'_{{A}}(a)=\delta_{{A}}(a)+r_0\otimes T(a)-T(a)\otimes r_0.
\end{eqnarray}
\end{Definition}

\begin{Definition}
Two flag datums $(\lambda, D, a_0, Q)$ and $(\lambda', D', a_0', Q') \in \mathcal{F}^{(2)}({A})$ are called equivalent if $\lambda'=\lambda$, $a_0'=a_0$ and there exist some element $r_0=r(x)\in {A}$ and $s\in k^{*}$ such that
\begin{eqnarray}
\label{eqq21}&&\lambda(a)r_0=[r_0,a] - D(a)+s D'(a),\\
\label{eqq22}&&\delta_{{A}}(r_0)+s Q'=r_0\otimes a_0-a_0\otimes r_0+Q.
\end{eqnarray}
\end{Definition}

From the above discussion, we obtain:
\begin{theorem}
Let $({A},[\cdot, \cdot], \delta_{{A}})$ be a Hom-Lie bialgebra of codimension one in a Hom-vector space $(E,\al_E)$. Then we have
$BExtd(E,{A})\cong \mathcal{HLB}(V,{A})\cong \mathcal{F}({A})/\equiv$.
\end{theorem}

\vskip7pt
\footnotesize{
\noindent Tao Zhang\\
College of Mathematics and Information Science,\\
Henan Normal University, Xinxiang 453007, P. R. China;\\
 E-mail address: \texttt{{zhangtao@htu.edu.cn}}

\end{document}